\documentclass[11pt]{article}
\usepackage{latexsym}
\usepackage{graphicx}
\usepackage{amssymb}

\newtheorem{theorem}{Theorem}[section]
\newtheorem{lemma}[theorem]{Lemma}
\newtheorem{proposition}[theorem]{Proposition}

\newtheorem{example}[theorem]{Example}

\newtheorem{remark}{Remark}[section]

\newenvironment{proof}{\smallskip\noindent{Proof:}}{\hfill\qed \newline}
\newcommand{\qed}{\quad\rule{1.5ex}{1.5ex}\vspace{1mm}\smallskip}
\newcommand{\forme}[1]{}

\newcommand{\text}{\mathsf}

\newcommand{\mc}{\mathcal}
\newcommand{\mbb}{\mathbb}

\title{Terwilliger Algebras of Wreath Powers\\ of One-Class Association Schemes}

\begin{document}

\author{Gargi Bhattacharyya and Sung Y. Song\\
Department of Mathematics, Iowa State University, Ames, Iowa, 50011,
USA}

\maketitle

\begin{abstract}
In this paper, we study the subconstituent algebras, also called as
Terwilliger algebras, of association schemes that are obtained as
the wreath product of one-class association schemes $K_n=H(1, n)$
for $n\ge 2$. We find that the $d$-class association scheme
$K_{n_{1}}\wr K_{n_{2}} \wr \cdots \wr K_{n_{d}}$ formed by taking
the wreath product of $K_{n_{i}}$ has the triple-regularity
property. We determine the dimension of the Terwilliger algebra for
the association scheme $K_{n_{1}}\wr K_{n_{2}}\wr \cdots \wr K_
{n_{d}}$. We give a description of the structure of the Terwilliger
algebra for the wreath power $(K_n)^{\wr d}$ for $n \geq 2$ by
studying its irreducible modules. In particular, we show that the
Terwilliger algebra of $(K_n)^{\wr d}$ is isomorphic to
$M_{d+1}(\mbb{C})\oplus M_1(\mbb{C})^{\oplus \frac12d(d+1)}$ for
$n\ge3$, and $M_{d+1}(\mbb{C})\oplus M_1(\mbb{C})^{\oplus
\frac12d(d-1)}$ for $n=2$.

\end{abstract}

%{\small {\it Keywords:} }

\section{Introduction}
The subconstituent algebra, which is also known as the Terwilliger
algebra of an association scheme was introduced by Terwilliger in
1992 as a new algebraic tool for the study of association schemes
\cite{Te92}. The Terwilliger algebra of a commutative association
scheme is a finite dimensional, semi-simple $\mbb{C}$-algebra, and
is noncommutative in general. This algebra helps understanding the
structure of the association schemes. It has been studied
extensively for many classes of association schemes. For example,
the Terwilliger algebra for $P$- and $Q$-polynomial association
schemes has been studied in \cite{Te93, TY94, Te96, Ca99}. The
structure of Terwilliger algebra of group association schemes has
been studied in \cite{BM95} and \cite{BO94}. In \cite{LMP06} the
structure of the Terwilliger algebra of a Hamming scheme $H(d, n)$
is given as symmetric $d$-tensors of the Terwilliger algebra of
$H(1,n)$ which are all isomorphic for $n>2$. It is also shown that
the Terwilliger algebra of $H(d,n)$ is decomposed as a direct sum of
Terwilliger algebra of hypercubes $H(d, 2)$ in \cite{LMP06}. They
deduce the decomposition into simple bilateral ideals using the
representation of classical groups. There is a detailed study of the
irreducible modules of the algebra for $H(d,2)$ in \cite{Go02}, and
for Doob schemes (the schemes coming as direct products of copies of
$H(2, 4)$ and/or Shrikhande graphs,) in \cite{Ta97}. Both of these
studies used elementary linear algebra and module theory.

In this paper, we study the Terwilliger algebras of association
schemes which are obtained as wreath products of $H(1, n)$, also
denoted $K_n$, for $n\ge 2$. We find that the $d$-class association
scheme $K_{n_{1}}\wr K_{n_{2}} \wr \cdots \wr K_{n_{d}}$ formed by
taking the wreath product of one-class association schemes
$K_{n_{i}}$ has the triple-regularity property in the sense of
\cite{Mu93} and \cite{Ja95}. Based on this fact, we determine the
dimension of the Terwilliger algebra for the association scheme
$K_{n_{1}}\wr K_{n_{2}}\wr \cdots \wr K_{n_{d}}$. We then find that
the wreath power $(K_n)^{\wr d}=K_n \wr K_n \wr \cdots \wr K_n$, $d$
copies of $K_n$, is formally self-dual and the Terwilliger algebra
is isomorphic to $M_{d+1}(\mbb{C})\oplus M_1(\mbb{C})^{\oplus
\frac12d(d+1)}$ for $n\ge3$, while $M_{d+1}(\mbb{C})\oplus
M_1(\mbb{C})^{\oplus \frac12d(d-1)}$ for $n=2$ in the notion of
Wedderburn-Artin's decomposition theorem of semisimple algebra. The
case $(K_2)^{\wr d}$ behaves a little differently from the general
case $(K_n)^{\wr d}$ for $n \geq 3$. We also study the corresponding
irreducible modules for $n=2$ in the course. We give the
decomposition by following the elementary approach employed in
\cite{TY94, Ta97, Go02}.

The remainder of the paper is organized as follows. In Section 2, we
provide the notation and terminology as well as a few basic facts on
the Terwilliger algebra and wreath product of association schemes
that will be used throughout. In Section 3, we discuss the structure
of the wreath product of one-class association schemes and compute
the dimension of its Terwilliger algebra. In Section 4, we study
Terwilliger algebras of wreath powers of one-class association
schemes and their irreducible modules. In Section 5, we make a
concluding remark and mention a few related open problems.

%\tableofcontents

\section{Preliminaries}

In this section we first briefly recall the notation and some basic
facts about association schemes and the Terwilliger algebra of a
scheme that are needed to deduce our results. Then we recall the
definition of wreath product of association schemes. For more
information on the topics covered in this section, we refer the
reader to \cite{BI84, BCN89, Te92, Te93, So02}.

\subsection{Association schemes and their Terwilliger algebras}
Let $X$ denote an $n$-element set, and let $M_X(\mbb{C})$ denote the
$\mbb{C}$-algebra of matrices whose rows and columns are indexed by
$X$. Let $\mathcal{X}=(X,\{R_{i}\}_{0\leq i \leq d})$ be a $d$-class
commutative association scheme of order $n$. The $d+1$ relations
$R_i\subseteq X\times X:=\{(x, y): x,y\in X\}$ are conveniently
described by their $\{0,1\}$-adjacency matrices $A_0, A_1,\dots,
A_d$ defined by $(A_i)_{xy}=1$ if $(x,y) \in R_i$; $0$ otherwise.
The \textit{intersection numbers} $p_{ij}^{h}$ are defined in terms
of the relations for the scheme by
$$p_{ij}^{h} =|\{ z \in X : (x,z) \in R_{i}, (z,y) \in R_{j}\}|$$
where $(x,y)$ is a member of the relation $R_h$. The definition of
an association scheme is equivalent to the following four axioms:
\begin{enumerate}
\item[(1)] $A_0= I$,
\item[(2)] $A_0 + A_1 + \cdots + A_d= J$,
\item[(3)] $A_i^{t}= A_{i'}$ for some $i^{'} \in \{ 0,1,\dots,d\}$,
\item[(4)] $A_iA_j = \sum\limits_{h=0}^{d} p_{ij}^{h}A_h$, for all
$i,j\in\{0, 1, \dots, d\}$
\end{enumerate}
where $I=I_n$ and $J=J_n$ are the $n \times n$ identity matrix and
all-ones matrix, respectively, and $A^{t}$ denotes the transpose of
the matrix $A$. The scheme is symmetric if $A_{i^{'}}=A_i$ for all
$i$, and is commutative if $p_{ij}^{h}=p_{ji}^{h}$ for all $h, i,
j$; and thus $A_iA_j=A_jA_i$.

%\subsection{The Bose-Mesner algebra}
If the scheme is commutative, the adjacency matrices generate a
$(d+1)$-dimensional commutative subalgebra $\mathcal{M}=\langle A_0
, A_1 , \dots , A_d \rangle$ of the full matrix algebra
$M_X(\mbb{C})$ over the field of complex numbers $\mbb{C}$. The
algebra $\mathcal{M}$ is known as the Bose-Mesner algebra of the
scheme. The Bose-Mesner algebra for a commutative association
scheme, being semi-simple, admits a second basis $E_0 , E_1 , \dots
, E_d$ of primitive idempotents. Note that $\mathcal{M}$ is also
closed under the Hadamard (entrywise) multiplication ``$\circ$" of
matrices. So there are nonnegative real numbers $q_{ij}^{h}$ called
the Krein parameters, such that
$$(4)\ E_i \circ E_j = |X|^{-1}\sum_{h=0}^d q_{ij}^{h}E_h.$$
There exist two sets of the $(d+1)^2$ complex numbers $p_j(i)$ and
$q_j(i)$ according to the $d+1$ expressions
$$A_j=\sum_{i=0}^{d}p_j(i)E_i,\quad \mbox{and} \quad E_j=|X|^{-1}\sum_{i=0}^{d}q_j(i)A_i.$$
The number $p_j(i)$ is characterized by the relation
$A_jE_i=p_j(i)E_i$. That is $p_j(i)$ is the eigenvalue of $A_j$,
associated with the eigenspace spanned by the columns of $E_i$,
occurring with the multiplicity $m_i=$ rank$(E_i)$. We define $P$ to
be the $(d+1) \times (d+1)$ matrix whose $(i,j$)-entry is $p_j(i)$.
This $P$ is referred to as the character table or first eigenmatrix
of the scheme $\mathcal{X}$.

Given an $n$-element set $X$, and the $\mbb{C}$-algebra
$M_X(\mbb{C})$ (or $M_n(\mbb{C})$), by the standard module of $X$,
we mean the $n$-dimensional vector space
$V=\mbb{C}^X=\bigoplus\limits_{x \in X}\mbb{C} \hat{x}$ of column
vectors whose coordinates are indexed by $X$. Here for each $x \in
X$, we denote by $\hat{x}$ the column vector with $1$ in the $x$th
position, and $0$ elsewhere. Observe that $M_X(\mbb{C})$ acts on $V$
by left multiplication. We endow $V$ with the Hermitian inner
product defined by $\langle u,v \rangle=u^{t}\overline{v}$ ($u,v \in
V$). For a given association scheme $\mathcal{X}=(X,\{R_{i}\}_{0\leq
i \leq d})$, $V$ can be written as the direct sum of $V_i=E_iV$
where $V_i$ are the maximal common eigenspaces of $A_0 , A_1 , \dots
, A_d$. Given an element $x \in X$, let $R_i(x)= \{y \in X: (x,y)
\in R_i\}$. The valencies $k_0,k_1,\dots,k_d$ of $\mathcal{X}$ are
denoted by $k_i=|R_i(x)|= p_{ii^{'}}^{0}.$ Let
$V_{i}^{*}=V_{i}^{*}(x) =\bigoplus\limits_{y \in R_i(x)}\mbb{C}
\hat{y}$. Both $R_i(x)$ and $V_i^*$ are referred to as the $i$th
\textit{subconstituent} of $\mathcal{X}$ with respect to $x$. Let
$E_i^{*}=E_i^{*}(x)$ be the orthogonal projection map from
$V=\bigoplus\limits_{i=0}^{d}V_{i}^{*}$ to the $i$th subconstituent
$V_{i}^{*}$. So, $E_i^{*}$ can be represented by the diagonal matrix
given by
$$(E_i^{*})_{yy}=\left
\{\begin{array}{lc}1 &\mbox{if } (x,y)\in R_i\\ 0 & \mbox{if } (x,y)
\notin  R_i \end{array}\right ..$$ The matrices
$E_0^{*},E_1^{*},\dots,E_d^{*}$ are linearly independent, and they
form a basis for a subalgebra $\mathcal{M}^{*}=\mathcal{M}^{*}(x) =
\langle E_0^{*},E_1^{*},\dots,E_d^{*} \rangle$ of $M_X(\mbb{C})$.

Let $A_i^{*}=A_i^{*}(x)$ be the diagonal matrix in $M_X(\mbb{C})$
with $yy$ entry $(A_i^{*})_{yy}=n\cdot (E_i)_{yx}$. Then we have
\begin{enumerate}
\item[(1)] $A_0^{*}=I$, \item[(2)] $A_0^{*}+A_1^{*}+\cdots+A_d^{*}=nE_0^{*}$,
\item[(3)]
$A_i^{*}A_j^{*}=\sum\limits_{h=0}^{d}q_{ij}^{h}A_h^{*}=A_jA_i$ for
all $i,j$.\end{enumerate} Furthermore, we also have
$$A_j^{*}=\sum_{i=0}^{d}q_j(i)E_i^{*} \quad \mbox{and}\quad
E_j^{*}=\frac1n\sum_{i=0}^{d}p_j(i)A_i^{*}.$$ Thus,
$A_0^{*},A_1^{*},\dots,A_d^{*}$ form a second basis for
$\mathcal{M}^{*}$. The algebra $\mathcal{M}^{*}$ is shown to be
commutative, semi-simple subalgebra of $M_X(\mbb{C})$. This algebra
is called the dual Bose-Mesner algebra of $\mathcal{X}$ with respect
to $x$.

%\subsection{Terwilliger algebra}

Let $\mc{X}=(X,\{R_i\}_{0 \leq i \leq d})$ denote a scheme and fix a
vertex $x \in X$, and let $\mathcal{T}=\mathcal{T}(x)$ denote the
subalgebra of $M_X(\mbb{C})$ generated by the Bose-Mesner algebra
$\mc{M}$ and the dual Bose-Mesner algebra $\mc{M}^{*}$. We call
$\mathcal{T}$ the {\it Terwilliger algebra of $\mc{X}$ with respect
to $x$}.
%By a {\it $\mc{T}$-module}, we mean a subspace $W \subseteq
%\mbb{C}^{X}$ such that $BW \subseteq W$ for all $B \in \mathcal{T}$.
%A $\mathcal{T}$-module is said to be irreducible whenever $W$ is
%nonzero and $W$ contains no other $\mc{T}$-modules other than $0$
%and $W$. By a $\mathcal{T}$-\textit{isomorphism} from $W$ to
%$W^{'}$, we mean a vector space isomorphism $\sigma : W \rightarrow
%W^{'}$ such that ($\sigma B - B \sigma)W=0$ for all $B \in
%\mathcal{T}$. The modules $W,W^{'}$ are said to be
%$\mathcal{T}$-\textit{isomorphic} whenever there exists a
%$\mathcal{T}$-isomorphism from $W$ to $W^{'}$. The standard module
%decomposes into irreducible $\mathcal{T}$-modules in a manner that
%reflects the structure of $\mathcal{T}$.
%By the \textit{center} of $\mathcal{T}$, we mean the subalgebra of
%$\mathcal{T}$ given by
%$${\rm{Center}}(\mathcal{T}):= \{C \in \mathcal{T} : BC =CB, \mbox{ for all }
%B \in \mathcal{T}\}.$$
The facts discussed in the rest of the
section will be useful when we describe the irreducible
$\mathcal{T}$-modules for our schemes in the subsequent sections.

%\subsection{The central primitive idempotents of $\mathcal{T}$}

\begin{lemma}\cite{Te92} Let $\mc{X}=(X,\{R_i\}_{0 \leq i \leq d})$
denote a $d$-class association scheme. For an arbitrary fixed vertex
$x \in X$, let $\mathcal{T}=\mathcal{T}(x)$. There exists a set
$\Phi = \Phi(x)$ and a basis $\{e_{\lambda}: \lambda \in \Phi \}$ of
the center of $\mc{T}$ such that
\begin{itemize}
\item[(i)] $I = \sum_{\lambda \in \Phi}e_{\lambda}$,

\item[(ii)] $e_{\lambda}e_{\mu} = \delta_{\lambda \mu}e_{\lambda}$
     (for all $\lambda, \mu \in \Phi$).\end{itemize}
\end{lemma}
We refer to $e_{\lambda}$ as the \textit{central primitive
idempotents} of $\mathcal{T}$. Let $V=\mbb{C}^X$ denote the standard
module.
\begin{lemma}\cite{Te92}\label{central} Let $\mc{X}=(X,\{R_i\}_{0 \leq i \leq d})$
denote a scheme and fix a vertex $x \in X$, and let
$\mathcal{T}=\mathcal{T}(x)$. Let $\{e_{\lambda}: \lambda \in \Phi
\}$ be the central primitive idempotents of $\mathcal{T}$.
\begin{enumerate}
\item[(i)] $V = \sum_{\lambda \in \Phi} e_{\lambda}V$ (Orthogonal
direct sum). Moreover, $e_{\lambda} : V \rightarrow e_{\lambda}V$ is
an orthogonal projection for all $\lambda \in \Phi$.
\item[(ii)] For each irreducible $\mathcal{T}$-module $W$, there is a
unique $\lambda \in \Phi$ such that $W \subseteq e_{\lambda}V$. We
refer to $\lambda$ as the \textit{type} of $W$.
\item[(iii)] Let $W$ and $W^{'}$ denote irreducible
$\mathcal{T}$-modules. Then $W$ and $W^{'}$ are
$\mathcal{T}$-isomorphic if and only if $W$ and $W^{'}$ have the
same type.
\item[(iv)] For all $\lambda \in \Phi,\ e_{\lambda}V$ can be decomposed
as an orthogonal direct sum of irreducible $\mathcal{T}$-modules of
type $\lambda$.
\item[(v)] Referring to (iv), the number of irreducible
$\mathcal{T}$-modules in the decomposition is independent of the
decomposition. We shall denote this number by mult$(e_{\lambda})$
(or simply mult$(\lambda)$) and refer to it as the multiplicity (in
$V$) of the irreducible $\mathcal{T}$-module of type $\lambda$.
\end{enumerate}
\end{lemma}

The set of triple products $E_{i}^{*}A_jE_{h}^{*}$ in $\mc{T}(x)$
plays a special role in our study, so we look at them little
closely. We can view $E_{i}^{*}A_jE_{h}^{*}$ as a linear map from
$V_h^{*} \rightarrow V_i^{*}$ such that
$$E_{i}^{*}A_jE_{h}^{*}\hat{y}=\sum_{z \in R_i(x) \cap
R_j(y)}\hat{z}$$ for each $\hat{y} \in V_h^{*}$. Terwilliger proved
the following key fact in \cite[Lemma 3.2]{Te92}.

\begin{proposition} \label{tripro} For $0 \leq h, i, j \leq d$,
$E_{i}^{*}A_jE_{h}^{*}=0$
if and only if $p_{ij}^{h}= 0$.\end{proposition}

Note that for every $i\in \{0, 1, \dots, d\}$, $A_i$ and $E_{i}^{*}$
can be written in terms of the triple products
$E_{i}^{*}A_jE_{h}^{*}$. Thus, the triple products
$E_{i}^{*}A_jE_{h}^{*}$ generate the Terwilliger algebra. It is
often easier to find the irreducible modules if we work with the
triple products $E_{i}^{*}A_jE_{h}^{*}$ instead of $A_i$ and
$E_{i}^{*}$.

\subsection{The Wreath product of association schemes}

We briefly recall the notion of the wreath product. Let
$\mathcal{X}=(X, \{R_i\}_{0\le i\le d})$ and $\mathcal{Y}=(Y,
\{S_j\}_{0\le j\le e})$ be association schemes of order $|X|= m$ and
$|Y|= n$. The wreath product $\mathcal{X} \wr \mathcal{Y}$ of
$\mathcal{X}$ and $\mathcal{Y}$ is defined on the set $X\times Y$;
but we take $Y = \{ y_1, y_2,\dots, y_n\}$, and regard $X\times Y$
as the disjoint union of $n$ copies $X_1,X_2,\dots, X_n$ of $X$,
where $X_j = X\times \{y_j\}$. The relations on $X_1 \cup X_2 \cup
\cdots \cup X_n$ is defined by the following rule:
\begin{quote} For any $j$, the relations between the elements of
    $X_j$ are determined by the association relations
    between the first coordinates in $\mathcal{X}$. For any $i$ and $j$,
the relations between $X_i$ and $X_j$ are determined by
    the association relation of the second coordinates
    $y_i$ and $y_j$ in $ \mathcal{Y}$ and the relation is
    independent from the first coordinates.\end{quote}

\noindent We may arrange the relations $W_0,W_1,\dots,W_{d+e}$ of
$\mathcal{X} \wr \mathcal{Y}$ as follows:
\begin{itemize}
\item $W_0 = \{ ((x,y),(x,y)) : (x,y)\in X \times Y\} $
\item $W_k = \{ ((x_1,y),(x_2,y)) : (x_1,x_2)\in R_k , y
    \in Y\} $ , for $1 \leq k \leq d$ ; and
\item $W_k = \{ ((x_1,y_1),(x_2,y_2)) : x_1,x_2 \in X,
    (y_1,y_2) \in S_{k-d}\} $ for $d+1 \leq k \leq d+e$.

\end{itemize}
\noindent Then the wreath product $\mathcal{X} \wr \mathcal{Y}= (X
\times Y, \{W_k\}_{0\le k\le d+e})$ is a $(d+e)$-class association
scheme. It is clear that $\mathcal{X} \wr \mathcal{Y}$ is
commutative (resp. symmetric) if and only if $\mathcal{X}$ and
$\mathcal{Y}$ are. Let $A_0,A_1,\dots, A_d$ and $C_0,C_1,\dots, C_e$
be the adjacency matrices of $\mathcal{X}$ and those of
$\mathcal{Y}$, respectively. Then the adjacency matrices $W_k$ of
$\mathcal{X} \wr \mathcal{Y}$ are given by
$$W_0=C_0 \otimes A_0, W_1=C_0\otimes A_1,\dots, W_{d}=C_0\otimes A_d,
W_{d+1}=C_1\otimes J_m, \dots, W_{d+e}=C_e\otimes J_m,$$ where
``$\otimes$" denotes the Kronecker product: $A\otimes B = (a_{ij}B)$
of two matrices $A=(a_{ij})$ and $B$. With the above ordering of the
association relations of $\mathcal{X} \wr \mathcal{Y}$, the relation
table of the wreath product is described by
 $$ R(\mathcal{X} \wr \mathcal{Y})=\sum\limits_{k=0}^{d+e}k\cdot
 A_k = I_n \otimes R(\mathcal{X}) + \{R(\mathcal{Y})+
 d(J_n-I_n)\}\otimes J_m.$$

\section{The Dimension of the Terwilliger Algebra of $K_{n_{1}}\wr
K_{n_{2}}\wr \cdots \wr K_ {n_{d}}$}

In this section, we first calculate the dimension of the Terwilliger
algebra of the wreath product of one-class association schemes. We
also show that these wreath product schemes satisfy the triple
regularity property.

Throughout the paper, for the notational simplicity, by $[n]$ we
denote the set of integers $\{1, 2, \dots, n\}$, and by $K_n=H(1,
n)$ we denote both the one-class association scheme $([n], \{R_0,
R_1\})$ with $A_1=J-I$, and the complete graph on $n$ vertices.

Let $\mc{X}=(X,\{R_i\}_{0\le i \le d})$ denote the $d$-class
association scheme $K_{n_{1}}\wr K_{n_{2}}\wr \cdots \wr K_{n_{d}}$
with
$$X=[n_1]\times [n_2]\times \cdots \times [n_d] \ =\ \{(a_1, a_2,
\dots, a_d): a_i\in [n_i], \mbox{ for } i=1, 2, \dots, d\}.$$ Let
$(1, 1, \dots, 1)\in X$ be a fixed base vertex $x$ of $\mc{X}$.
Without loss of generality, we can arrange the association relations
such that
\begin{itemize}
%\item $R_0(x)=\{x\}$,
\item $R_1(x)=\{(a, 1, 1, \dots, 1): a\in \{2,3,\dots, n_1\}\}$,
\item for $i=2, 3, \dots, d$, $$R_i(x)=\{(a_1, a_2, \dots, a_{i-1}, b, 1, 1, \dots, 1):
a_k\in [n_k] \mbox{ for } k\in \{1, 2, \dots, i-1\}, b\in [n_i]-
\{1\}\}.$$
\end{itemize}
\noindent We observe that $k_i=|R_i(x)|=(n_i
-1)\prod\limits_{k=1}^{i-1}n_k$. We can arrange rows and columns of
the relation table of $\mc{X}$ by the order of parts in the
partition $X=R_0(x)\cup R_1(x)\cup \cdots \cup R_d(x)$.

\begin{example}
The following is the relation table for the wreath product of three
association schemes $K_2$, $K_2$ and $K_3$ arranging the elements in
the order described above.
\[R(K_2 \wr K_2 \wr K_3)= \left(\begin{array}{cc|cc|cccccccc}

0 & 1& 2& 2& 3& 3& 3& 3& 3& 3& 3& 3\\
1& 0& 2 & 2 & 3& 3& 3& 3& 3& 3& 3& 3\\ \hline
2 & 2&  0 & 1& 3& 3& 3& 3& 3& 3& 3& 3\\
2 & 2 &  1& 0& 3& 3& 3& 3& 3& 3& 3& 3\\ \hline

3 & 3 & 3 & 3 & 0 & 1& 2& 2&  3& 3& 3& 3\\
3 & 3 & 3 & 3 & 1& 0& 2 & 2 & 3& 3& 3& 3\\
3 & 3 & 3 & 3 & 2 & 2&  0 & 1&  3& 3& 3& 3\\
3 & 3 & 3 & 3 & 2 & 2 &  1& 0&  3& 3& 3& 3\\

3& 3& 3& 3& 3& 3& 3& 3&0 & 1& 2& 2\\
3& 3& 3& 3& 3& 3& 3& 3&1& 0& 2 & 2 \\
 3& 3& 3& 3& 3& 3& 3& 3&2 & 2&  0 & 1\\
 3& 3& 3& 3& 3& 3& 3& 3&2 & 2 &  1& 0\\
\end{array}\right)\]
We note that any table obtained from this table by permuting the
order of rows and corresponding columns simultaneously represents
the same association scheme.
\end{example}

\begin{lemma}\label{non-0-para}
Let $\mathcal{X}=(X, \{R_i\}_{0\le i\le d})= K_{n_{1}}\wr
K_{n_{2}}\wr \cdots \wr K_{n_{d}}$.  Then the complete list of
nonzero $p^h_{ij}$ where $h, i, j \in \{0,1,2, \dots ,d\}$ is as
follows.
\begin{enumerate}
\item[$(1)$] For $h=0$, $k_0=p_{00}^0=1,\ k_1=p_{11}^0=n_1-1$, and\\
$k_j=p_{jj}^0=(n_j-1)\prod\limits_{l=1}^{j-1}n_l$ \quad for $j=2, 3,
\dots, d$.
\item[$(2)$] For $h=1, 2, \dots, d$,
\begin{enumerate}
\item[$(a)$] $p^h_{hh}=(n_h-2)\prod\limits_{l=1}^{h-1}n_l$.
\item[$(b)$] $p^h_{jj}=(n_j-1)\prod\limits_{l=1}^{j-1}n_l$\quad for $h+1 \leq j\le d$
\item[$(c)$] $p^h_{jh}=p^h_{hj}=(n_j-1)\prod\limits_{l=1}^{j-1}n_l$\quad for $1\le j \leq h-1$
\item[$(d)$] $p^h_{0h}=p^h_{h0}=1$
\end{enumerate}
\end{enumerate}
\end{lemma}
\begin{proof} It is straightforward to calculate the intersection
numbers.
\end{proof}
Due to this lemma, we have the following list of non-zero triple
products in $\mc{T}$.
\begin{theorem} \label{non-0-tri-pro} The complete list of nonzero
triple products $E_i^*A_jE_h^*$ among all $h, i, j \in \{0,1,2
\cdots ,d\}$ in $\mathcal{X}= K_{n_{1}}\wr K_{n_{2}}\wr \cdots \wr
K_{n_{d}}$ is given as follows.
\begin{enumerate}
\item[$(1)$] $E_i^*A_iE_0^*$ for $0\le i\le d$
\item[$(2)$] $E_h^*A_hE_h^*$ if and only if $n_h\ge 3$ for $1 \le h \leq
d$,
\item[$(3)$] $E_j^*A_jE_h^*$ for $2\le h+1 \leq j\le d$
\item[$(4)$] $E_j^*A_hE_h^*$ for $1\le j+1 \leq h \leq d$
\item[$(5)$] $E_h^*A_jE_h^*$ for $1\le j+1\le h\le d$
\end{enumerate}
\end{theorem}
\begin{proof} Immediate from the above lemma by Proposition
\ref{tripro}.
\end{proof}

In order to calculate the dimension of the Terwilliger algebra for
the $d$-class association scheme $\mathcal{X}=(X, \{R_i\}_{0\le i\le
d})= K_{n_{1}}\wr K_{n_{2}}\wr \cdots \wr K_{n_{d}}$. Let $x\in X$
be a fixed vertex, and consider the subspace
$\mc{T}_0=\mathcal{T}_0(x)$ of $\mc{T}=\mathcal{T}(x)$ spanned by
$\{E_i^*A_jE_h^*: 0 \leq i,j,h \leq d\}$. It is easy to see that
$\mathcal{T}$ is generated by $\mathcal{T}_0$ as an algebra since
$\mathcal{T}_0$ contains $A_i$ and $E_i^{*}$ for all $i$, but in
general, $\mathcal{T}_0$ may be a proper linear subspace of
$\mathcal{T}$. However, we will see that $\mc{T}_0 =
 \mc{T}$ for $K_{n_{1}}\wr K_{n_{2}}\wr \cdots \wr K_{n_{d}}$ shortly.
First we have the following dimension formula for $\mc{T}_0$.

\begin{theorem}\label{dim} Let $\mathcal{X}=
K_{n_{1}}\wr K_{n_{2}}\wr \cdots \wr K_{n_{d}}$. Then the dimension
of $\mathcal{T}_0$ is given by $$\dim(\mc{T}_0)=(d+1)^2+\frac12
d(d+1) - b$$ where $b$ is the number of $K_2$ factors in the wreath
product. In particular,
$$ (d+1)^2+\left (\begin{array}{c} d\\ 2\\ \end{array}\right ) \leq
dim(\mathcal{T}_0) \leq (d+1)^2 +\left (\begin{array}{c} d+1\\ 2\\
\end{array}\right )$$

\end{theorem}

\begin{proof} In the Theorem \ref{non-0-tri-pro}, the number of non-zero
triple products can be counted as d+1 from (1), $d-b$ from (2),
$\frac12 d(d-1)$ from (3), and $d(d+1)$ from (4) and (5). As they
are independent of each other we have the $dim(\mathcal{T}_0(x))=d
+1 + d-b+ \frac12 d(d-1)+d(d+1)=(d+1)^2+\frac12 d(d+1)-b$ as
desired. The case when all $n_i=2$ gives the lower bound as in this
case $b=d$. The upper bound is given by the case where $n_i\ge 3$
for all $i$. In such a situation $b=0$. This completes the proof.
\end{proof}

We now show that $\mc{T}=\mc{T}_0$ for the wreath product scheme
$K_{n_{1}}\wr K_{n_{2}}\wr \cdots \wr K_{n_{d}}$; so the scheme has
the triple-regularity property. The concept of triple-regularity was
first studied by Terwilliger. For more information on it, we refer
to \cite[p.120]{Ja95}). We use the following equivalent properties
of triple-regularity observed by Munemasa.

\begin{proposition}[\cite{Mu93}]\label{tri-reg} Let $\mathcal{X}$ be a commutative association
scheme. Then the following are equivalent. \begin{enumerate}
\item[$1.$] $\mathcal{X}$ is triply regular; i.e. $\mc{X}$ has the property that
the size of the set $R_i(x)\cap R_j(y)\cap R_h(z)$ depends only on
the set $\{i,j,h,l,m,n\}$ where $(x,y) \in R_l$, $(x,z) \in R_m$ and
$(y,z) \in R_n$.
\item[$2.$] $A_iE_j^{*}A_h \in \mathcal{T}_0$ for any $h, i, j$.
\item[$3.$] $\mc{T}(x)=\mc{T}_0(x)$ for $x\in X$.\end{enumerate}
\end{proposition}

According to this proposition, it suffices to verify that all triple
products $A_iE_h^*A_j$ belong to $\mc{T}_0$ in order to show that
$\mc{T}=\mc{T}_0$ for $K_{n_{1}}\wr K_{n_{2}}\wr \cdots \wr
K_{n_{d}}$.

\begin{lemma} \label{conv1} For the $d$-class scheme
$K_{n_{1}}\wr K_{n_{2}}\wr \cdots \wr K_{n_{d}}$, we have the
following.
\begin{enumerate}
\item[$(1)$] $A_iE_h^*A_j=(A_iE_h^*)(E_h^*A_j)$ for all $h, i, j\in
\{0, 1, \dots, d\}$.
\item[$(2)$] $A_hE_h^* =\sum_{j=0}^{h}E_j^{*}A_hE_h^{*}$
for all $h\in \{0, 1, \dots, d\}$.
\item[$(3)$] For $0\le i < h\le d$,
$A_iE_h^* = E_h^{*}A_iE_h^{*}$.
\item[$(4)$] For $0\le h< i\le d$, $A_iE_h^* = E_i^{*}A_iE_h^{*}$.
\end{enumerate}
\end{lemma}
\begin{proof}
(1) It is trivially true as $E_h^*$ are idempotents.\\
(2) The nonzero entries of $A_hE_h^*$ are the nonzero entries of the
columns of $A_h$ indexed by vertices in $R_h(x)$. The rest of the
entries are zero. The columns of $A_h$ indexed by vertices in
$R_h(x)$ have $1$ in rows indexed by the vertices in $R_{0}(x) \cup
R_{1}(x) \cup \cdots \cup R_{h-1}(x)$. In addition the rows indexed
by the vertices $R_{h}(x)$ have zero in the diagonal blocks of size
$(\prod_{i=1}^{h-1}n_i) \times (\prod_{i=1}^{h-1}n_i)$, and $1$
elsewhere. These are essentially
$\sum\limits_{j=0}^{h}E_j^{*}A_hE_h^{*}$.\\
(3) If $i < h$, the nonzero entries of $A_iE_h^*$ are the nonzero
entries of the columns of $A_i$ indexed
 by vertices in $R_h(x)$. The rest of the entries are zero. The columns
 of $A_i$ indexed
 by vertices in $R_h(x)$ have $1$ in rows indexed by the vertices in $R_h(x)$.
 The rest of the entries are zero.\\
(4) If $i > h$, the nonzero entries of $A_iE_h^*$ are the nonzero
entries of the columns of $A_i$ indexed
 by vertices in $R_h(x)$. The rest of the entries are zero. The columns
 of $A_i$ indexed
 by vertices in $R_h(x)$ have $1$ in rows indexed by the vertices in $R_i(x)$.
 The rest of the entries are zero. This completes the proof.
\end{proof}

\begin{lemma} \label{conv2} In $K_{n_{1}}\wr K_{n_{2}}\wr \cdots \wr K_{n_{d}}$,
\begin{enumerate}
\item[$(1)$] $E_h^*A_h=\sum_{j=0}^{h}E_h^{*}A_hE_j^{*}$ for $h\in \{0, 1, \dots, d\}$;
\item[$(2)$] for $0\le i < h\le d$, $E_h^*A_i = E_h^{*}A_iE_h^{*}$;
\item[$(3)$] for $0\le h< i\le d$, $E_h^*A_i =
E_h^{*}A_iE_i^{*}$.\end{enumerate}
\end{lemma}
\begin{proof} It follows from the previous lemma and the fact that the transpose
of $A_iE_h^*$ is $E_h^*A_i $.
\end{proof}

\begin{lemma}\label{conv-lc1}
For $K_{n_{1}}\wr K_{n_{2}}\wr \cdots \wr K_{n_{d}}$, we have the
following linear combinations for $A_iE_h^*A_j$.
\begin{enumerate}
\item[$(1)$] For $h \in\{2,3, \dots , d\}$,
\[A_hE_{h}^{*}A_h =
(n_h-1)\left(\prod\limits_{k=1}^{h-1}n_k\right)\sum_{m =
0}^{d}\sum_{l = 0}^{h-1}\sum_{n = 0}^{h-1}E_{l}^*A_{m}E_{n}^*\]
\[+\ (n_h-2)\left(\prod\limits_{k=1}^{h-1} n_k\right)\left
\{\sum_{m=0}^{d}\sum_{n=0}^{h-1}E_{h}^*A_{m}E_{n}^* +
\sum_{m=0}^{d}\sum_{l=0}^{h}E_{l}^*A_{m}E_{h}^*\right \}.\]
\item[$(2)$]
\begin{enumerate}
\item[$(a)$] For $d\ge j > h\ge 2$, \[A_{h}E_{h}^{*}A_{j}
= (n_h-1)\left(\prod\limits_{k=1}^{h-1}n_k\right)\sum_{m =
0}^{d}\sum_{l = 0}^{h-1}E_{l}^*A_{m}E_{j}^*\ +\
(n_{h}-2)\left(\prod\limits_{k=1}^{h-1}n_k\right)\sum_{m=0}^{d}E_{h}^*A_{m}E_{j}^*;\]
\item[$(b)$] for $2\le j< h\le d$,
$$A_{h}E_{h}^{*}A_{j} = (n_{j}-1)\left(\prod\limits_{k=1}^{j-1}n_k\right)\sum_{m =
0}^{d}\sum_{l = 0}^{h-1}E_{l}^*A_{m}E_{h}^*.$$
\end{enumerate}
\item[$(3)$]
\begin{enumerate} \item[$(a)$] For $d\ge i > h\ge 2$,
$$A_{i}E_{h}^{*}A_{h} =(n_{h}-1)\left(\prod\limits_{k=1}^{h-1}n_k\right)
\sum_{m = 0}^{d}\sum_{n = 0}^{h-1}E_{i}^*A_{m}E_{n}^*\ + \
(n_{h}-2)\left(\prod\limits_{k=1}^{h-1}n_k\right)\sum_{m=0}^{d}E_{i}^*A_{m}E_{h}^*;$$
\item[$(b)$] for $i < h$,
$$A_{i}E_{h}^{*}A_{h} = n_{1} \cdots n_{i-1}(n_{i}-1)\sum_{m =
0}^{d}\sum_{n = 0}^{h-1}E_{h}^*A_{m}E_{n}^*.$$
\end{enumerate}
\item[$(4)$]
\begin{enumerate}
\item[$(a)$] For $d\ge i > h\ge 2$, $$A_{i}E_{h}^{*}A_{i} = n_{1} \cdots
n_{h-1}(n_{h}-1)\sum_{m = 0}^{d}E_{i}^*A_{m}E_{i}^*;$$
\item[$(b)$] for $2\le i < h\le d$,
$$A_{i}E_{h}^{*}A_{i} =n_{1} \cdots
n_{i-1}(n_{i}-1)\sum_{m = 0}^{d}E_{h}^*A_{m}E_{h}^*.$$
\end{enumerate}
\end{enumerate}
\end{lemma}
\begin{proof} (1) Applying the identities in Lemmas \ref{conv1} and
\ref{conv2} to
 $A_{i}E_{h}^{*}A_{j}=(A_iE_h^*)(E_h^*A_j)$  we have
$$A_{i}E_{h}^{*}A_{j} = k_h \sum_{m = 0}^{d}\sum_{l = 0}^{h-1}\sum_{n =
0}^{h-1}E_{l}^*A_{m}E_{n}^*\qquad\qquad\qquad\qquad\qquad\qquad$$
$$\qquad\qquad + (k_h-(k_0+ \cdots +k_{h-1}))\left \{\sum_{m=0}^{d}\sum_{n=0}^{h-1}E_{h}^*A_{m}E_{n}^*
+ \sum_{m=0}^{D}\sum_{l=0}^{h}E_{l}^*A_{m}E_{h}^*\right \}$$ with
$k_h=n_1n_2\cdots n_h -n_1n_2\cdots n_{h-1}$ as desired.\\
(2) The proof of part (a) follows from Lemma \ref{conv1}(1) and
\ref{conv2}(3), while (b) follows from Lemma \ref{conv1}(1) and
\ref{conv2}(2).\\
(3) The proof is a similar to part (2).\\
(4) The proof of (a) follows from  Lemma \ref{conv1}(3) and
\ref{conv2}(3), and (b) follows from Lemma \ref{conv1}(2) and
\ref{conv2}(2).
\end{proof}

\begin{lemma}\label{conv-lc2}
In $K_{n_{1}}\wr K_{n_{2}}\wr \cdots \wr K_{n_{d}}$, for $i , j , h
\in\{0,1, \dots , d\}$, suppose that no two of $i, j, h$ are equal.
\begin{enumerate}
\item[$(i)$] If $i > h$ and $j > h$, then
 $$A_{i}E_{h}^{*}A_{j} = (n_{h}-1)\left(\prod\limits_{k=1}^{h-1}n_k\right)\sum_{m =
0}^{d}E_{i}^*A_{m}E_{j}^*.$$
\item[$(ii)$] If $i < h< j$, then
 $$A_{i}E_{h}^{*}A_{j} = (n_{i}-1)\left(\prod\limits_{k=1}^{i-1}n_k\right)\sum_{m = 0}^{d}E_{h}^*A_{m}E_{j}^*.$$
\item[$(iii)$] If $i > h > j$, then
 $$A_{i}E_{h}^{*}A_{j} = (n_{j}-1)\left(\prod\limits_{k=1}^{j-1}n_k\right)\sum_{m = 0}^{d}E_{i}^*A_{m}E_{h}^*.$$
\item[$(iv)$] If $i < h$ and $j < h$, then
\begin{enumerate} \item[$(a)$] for $i < j$,
 $$A_{i}E_{h}^{*}A_{j} = (n_{i}-1)\left(\prod\limits_{k=1}^{i-1}n_k\right)E_{h}^*A_{j}E_{h}^*;$$
\item[$(b)$] for $i > j$,
 $$A_{i}E_{h}^{*}A_{j} = (n_{j}-1)\left(\prod\limits_{k=1}^{j-1}n_k\right)E_{h}^*A_{i}E_{h}^*.$$
\end{enumerate}
\end{enumerate}
\end{lemma}

\begin{proof}(i),(ii) and (iii) are similar to Lemma 3.2.8

(iv)(a) Lemmas \ref{conv1}(ii) and \ref{conv2}(ii) gives us nonzero
entries of $A_{i}E_{h}^{*}A_{j}$ occur in the rows and columns
indexed by the vertices $R_h(x)$. Consider the diagonal blocks of
size $k_j \times k_j$ inside $A_{i}E_{h}^{*}$ indexed by the rows
and columns of vertices in $R_h(x)$. These diagonal blocks have
$k_i$ $1$'s in each row and column. The off diagonal entries are all
zero. In a similar manner consider diagonal blocks of size $k_j
\times k_j$ inside $A_{i}E_{h}^{*}A_{j}$ indexed by the rows and
columns of vertices in $R_h(x)$. These diagonal blocks are all zero
and the off diagonal entries are $1$. This observation gives us
$$A_{i}E_{h}^{*}A_{j} = n_{1} \cdots
n_{i-1}(n_{i}-1)E_{h}^*A_{j}E_{h}^*$$

(b) Similar to part (a).
\end{proof}

\begin{theorem} \label{trireg} The $d$-class association scheme $K_{n_{1}}\wr K_{n_{2}}\wr\cdots
\wr K_{n_{d}}$ is triply regular and $\mathcal{T}= \mathcal{T}_0$.
\end{theorem}

\begin{proof} It is straightforward to check that each of the rest of
$A_iE_h^*A_j$ that are not covered in Lemma \ref{conv-lc1} and Lemma
\ref{conv-lc2}, also can be expressed as linear combination of the
generators of $\mc{T}_0$. Thus the conclusion follows from
Proposition \ref{tri-reg}.
\end{proof}
In summary, we have the following.

\begin{theorem}\label{dimension} The dimension of the Terwilliger algebra
$\mathcal{T}$ of $K_{n_{1}}\wr K_{n_{2}}\wr\cdots\wr K_ {n_{d}}$ is
$$dim(\mc{T})=(d+1)^2 + \frac12 d(d+1)-b$$ where $b$ denotes the number of
factors $K_{n_i}$ with $n_i=2$.
\end{theorem}

\begin{proof} Immediate consequence of Lemma \ref{dim}, Proposition \ref{tri-reg},
and Theorem \ref{trireg}.
\end{proof}

\section{Terwilliger Algebras of Wreath Powers of $K_m$}

In this section we describe the structure of the Terwilliger algebra
of the wreath power $(K_m)^{\wr d}=K_m\wr K_m\wr \cdots \wr K_m$ of
$d$ copies of the one-class association scheme $K_m$. We begin the
section with the description of the Terwilliger algebra of the
one-class association scheme $K_m$ of order $m$. The nontrivial
relation graph of $K_m$ is the adjacency matrix of the complete
graph $K_m$ which is also viewed as Hamming graph $H(1,m)$. We then
describe the Terwilliger algebra of $(K_m)^{\wr 2}$ whose first
relation graph is the complete $m$-partite strongly regular graph
with parameters $(v, k, \lambda, \mu)=(m^2, m(m-1), m(m-2),
m(m-1))$. We then compare the combinatorial structures of wreath
square $(K_m)^{\wr 2}$ and cube $(K_m)^{\wr 3}$ of $K_m$ to describe
irreducible $\mc{T}$-modules of $(K_m)^{\wr 3}$ by extending those
of $(K_m)^{\wr 2}$. Similarly, the structure of the irreducible
$\mc{T}$-modules of $(K_m)^{\wr d}$ for any higher $d$ will be
described from that of $(K_m)^{\wr (d-1)}$. It is shown that all
non-primary irreducible $\mc{T}$-modules of wreath powers of $K_m$
are of dimension 1. We conclude the section by describing the
Terwilliger algebra of the $d$-power $(K_m)^{\wr d}$ for an
arbitrary $d\ge 2$.

\subsection{The Terwilliger algebra of $K_m$}
Let $K_m=([m], \{R_0, R_1\})$, and let $x=1$. Then
\[R_1(x)=\{2, 3, \dots, m\}\]
and, we may denote $A_1$ by
\[A_1=J-I=\left [\begin{array}{cc} 0 & \mathbf{1}^t\\
& \\
 \mathbf{1} & J_{m-1}-I_{m-1}\\ \end{array}\right ]\] where
$\mathbf{1}$ is the $(m-1)$-dimensional column vector all of whose
entries are 1.

\begin{remark} By Theorem \ref{dimension}, we know that the dimension
of the Terwilliger algebra of $K_m$ is $5$ if $m>2$ and $4$ if
$m=2$. Also by Theorem \ref{non-0-tri-pro}, all the matrices in the
Terwilliger algebra of $K_m$ is a linear combination of the matrices
$E_0^*A_0E_0^*$, $E_0^*A_1E_1^*$, $E_1^*A_1E_0^*$, $E_1^*A_0E_1^*$,
and $E_1^{*}A_1E_1^{*}$. (If $m=2$, then $E_1^{*}A_1E_1^{*}=0$.)
\end{remark}

If we set \[E_{11}=E_0^*A_0E_0^*,\ E_{12}=E_0^*A_1E_1^*,\
E_{21}=\frac{1}{m-1}E_1^*A_1E_0^*,\
E_{22}=\frac{1}{m-1}(E_1^*A_0E_1^*+E_1^{*}A_1E_1^{*}),\] then these
matrices form a subalgebra $\mc{U}$ of $\mc{T}(x)$ as its
multiplication table is given by
\[\begin{array}{c|cccc}
  & E_{11} & E_{12} & E_{21} & E_{22} \\ \hline
E_{11}& E_{11} & E_{12} & 0 & 0\\
E_{12} & 0 & 0 & E_{11} & E_{12}\\
E_{21} &E_{21} & E_{22} & 0 & 0\\
E_{22} & 0 & 0 & E_{21} & E_{22}\\
\end{array}\]
Considering the isomorphism between $\mc{U}$ and $M_2(\mathbb{C})$
that takes
\[E_{11}\mapsto \left[\begin{array}{cc}
1& 0\\
0& 0\\
\end{array}\right];\quad E_{12} \mapsto
\left[\begin{array}{cc}
0& 1\\
0& 0\\
\end{array}\right];\quad E_{21} \mapsto
\left[\begin{array}{cc}
0& 0\\
1& 0\\
\end{array}\right]; \quad E_{22} \mapsto
\left[\begin{array}{cc}
0& 0\\
0& 1\\
\end{array}\right].\]
we see that $\mc{T}(x)=M_2(\mbb{C})\oplus M_1(\mbb{C})$ by
Wedderburn-Artin's Theorem (cf. \cite[Sec. 2.4]{DK94}).
Specifically, if we set $F=E_1^{*}A_0E_1^*-E_{22}$, then it turns
out that $FX=0$ for all $X \in \mathcal{U}$. This gives us
$\mathcal{T}(x) = \mathbb{C}F \oplus \mathcal{U}$. While $F=0$ for
$m=2$, $F \neq 0$ for all $m>2$. Therefore, we reasserted the
following.
\begin{theorem}\cite{LMP06} The Terwilliger algebra of $K_m$ can be
described as follows:
\[
\mathcal{T}(x)\cong \left\{\begin{array}{ll}
M_2(\mathbb{C}) & \textrm{if  } m=2 \\
M_2(\mathbb{C}) \oplus M_1(\mathbb{C}) & \textrm{if  }m >2 \\
\end{array} \right..
\]
\end{theorem}

\subsection{The Terwilliger algebra of $(K_m)^{\wr 2}$ and its irreducible modules}

Let $\mathcal{X}= (K_m)^{\wr 2}=(X, \{R_0, R_1, R_2\})$ be the
wreath square of $K_m$. Without loss of generality, we arrange the
relations so that the first relation graph $(X, R_1)$ is to be the
complete $m$-partite strongly regular graph with parameters $(v,k,
\lambda, \mu)=(m^2, m(m-1), m(m-2), m(m-1))$, which is also the
wreath square of the complete graph $K_m$.

Let $X=[m]\times [m]=\{(i,j):\ i, j\in [m] \}$, and let $x=(1,1)$.
We will refer to $(1,1)$ as the base vertex $x$. Then %$R_i(x) = \{y
%\in X :\ (x,y) \in R_i\}$ are given by
\[R_1(x)=\{(i,j):\ i\in [m],\ j\in \{2, 3, \dots, m\}\},\]
\[R_2(x)=\{ (i, 1):\ i\in \{2, 3, \dots, m\} \}.\]
Let the adjacency matrices $A_i$ and the relation table $R$ of
$\mathcal{X}$ be decomposed according to the partition $X= R_0(x)
\cup R_1(x) \cup R_2(x)$. Then,
\[A_{1}= \left[\begin{array}{ccc}
0&\mathbf{1}_2^t& \mathbf{0}_1^t\\
\mathbf{1}_2 & B_{2} & L\\
\mathbf{0}_1 & L^t & B_{1}\\
\end{array}\right],\quad
A_{2}= \left[\begin{array}{ccc}
0&\mathbf{0}_2^t& \mathbf{1}_1^t\\
\mathbf{0}_2 & C_{2} & N\\
\mathbf{1}_1 & N^t & C_{1}\\
\end{array}\right],\quad
R= \left[\begin{array}{ccc}
0&\mathbf{1}_2^t& 2\mathbf{1}_1^t\\
\mathbf{1}_2 & B_2+2C_{2} & L\\
2\mathbf{1}_1 & L^t & 2C_{1}\\
\end{array}\right].\]
where $\mathbf{1}_2$ and $\mathbf{1}_1$ are all-ones column vectors
of size $m(m-1)$ and $(m-1)$, respectively; $\mathbf{0}_2$ and
$\mathbf{0}_1$ are all-zeros column vectors of size $m(m-1)$ and
$m-1$, respectively; $L$ and $N$ are $m(m-1) \times (m-1)$ all-ones
and all-zeros matrices, respectively;
$B_2=J_{m(m-1)}-(I_{m-1}\otimes J_{m})$, and $B_1$ is a $(m-1)
\times (m-1)$ zero matrix, while $C_2=I_{m-1}\otimes (J_m-I_m)$ and
$C_1=J_{m-1}-I_{m-1}$.

We note that $A_2 = \frac{1}{m(m-1)} A_1^2 - \frac{m-2}{m-1}A_1 -I$.
We also note that $(K_{m})^{\wr 2}$ is formally self-dual $P$- and
$Q$-polynomial association scheme with its first and second
eigenmatrices
\[P = Q= \left[\begin{array}{ccc}
1& m(m-1) & m-1\\
1 & 0 & -1\\
1 &-m & m-1\\
\end{array}\right].\]
The characteristic polynomial of $A_1$ is $\theta^2+(\mu
-\lambda)\theta +(\mu-k)=0$, and the eigenvalues of $A_1$ are given
by $k=m(m-1),\ r=0,\mbox{ and } s=-m$ with multiplicities, $1,\
m(m-1), \mbox{ and } m-1$, respectively. The induced subgraphs of
the graph $(K_m)^{\wr 2}$ on the vertex sets $R_1(x)$ and $R_2(x)$
with adjacency matrices $B_2$ and $B_1$ respectively, are called the
subconstituents of the graph with respect to $x$. Set
$\lambda=p_{11}^{1}=m(m-2)$, $\mu = p_{11}^{2}=m(m-1)$,
$k_1=m_1=p_{11}^{0}=m(m-1), \mbox{ and } k_2=m_2=m-1$. Furthermore,
$B_2\mathbf{1}_2=m(m-1)\mathbf{1}_2$ and
$C_1\mathbf{1}_1=(m-1)\mathbf{1}_1$. The eigenvalues of $B_2$ are
$m(m-2)$, $-m$ and $0$ with multiplicities $1$, $m-2$ and $m^2-2m+1$
respectively. For $B_1$, $0$ is the only eigenvalue with
multiplicity $m-1$.

Cameron, Goethals and Seidel introduced the concept of restricted
eigenvalues and eigenvectors \cite{CGS78}. An eigenvalue of $B_2$
(resp. $B_1$) is called restricted if it has an eigenvector
orthogonal to the all-ones vector of size $k_1$ (resp. $k_2$).
Tomiyama and Yamazaki used restricted eigenvectors, the eigenvectors
associated with the restricted eigenvalue of $B_i$ that are
orthogonal to $\mathbf{1}_i$, to describe the subconstituent algebra
of a 2-class association scheme constructed from a strongly regular
graph.  In order to describe the irreducible $\mc{T}$-modules of our
scheme $(K_{m})^{\wr 2}$, we will need the result of Theorem 5.1 in
\cite{CGS78} for our scheme.
\begin{lemma} With the above notations for $(K_{m})^{\wr 2}$ and
its eigenvalues $m(m-1)$, $0$ and $-m$, we have
the following.
\begin{enumerate}\label{res-eigen}
\item Suppose $\mathbf{y}$ is a restricted eigenvector of $B_2$
with an eigenvalue $\theta$. Then $\mathbf{y}$ is the eigenvector of
$LL^t$ with the eigenvalue $-\theta (\theta +m)$, and
$L^t\mathbf{y}$ is the zero vector or the restricted eigenvector of
$B_1$ with the eigenvalue $-m-\theta$. In particular,
$L^t\mathbf{y}$ is zero if and only if $\theta \in \{ 0,-m\}$.

\item Suppose $\mathbf{z}$ is a restricted eigenvector of $B_1$
with an eigenvalue $\theta'$. Then $\mathbf{z}$ is the eigenvector
of $L^tL$ with the eigenvalue $-\theta' (\theta' +m)$, and
$L\mathbf{z}$ is the zero vector or the restricted eigenvector of
$B_2$ with the eigenvalue $-m-\theta'$. In particular, $L\mathbf{z}$
is zero if and only if $\theta' \in \{ 0,-m\}$.
\end{enumerate}
\end{lemma}

\begin{proof} It is immediate from the fact that $B_2$ is the
adjacency matrix of the strongly regular graph with parameters
$(m(m-1),m(m-2),m(m-3),m(m-2))$.
\end{proof}
The next three lemmas are results of Tomiyama and Yamazaki (reported
in \cite{TY94}) suited to our association scheme $(K_{m})^{\wr 2}$.
\begin{lemma} Let $\mathcal{T}$ denote the Terwilliger algebra of
$(K_{m})^{\wr 2}$. Then with the above notations in this subsection,
we have the following.
\begin{enumerate}
\item[$(1)$] Suppose $\mathbf{y}$ is a restricted eigenvector of $B_2$ with an
eigenvalue $\theta$. Then the vector space $W$ over $\mbb{C}$ which
is spanned by $(0,\mathbf{y}^t,\mathbf{0}_1^t)^t$ is a thin
irreducible $\mathcal{T}$-module over $\mbb{C}$ and dim $W=1$ if
$\theta \in \{ 0,-m\}$.
\item[$(2)$] Suppose $\mathbf{z}$ is a restricted eigenvector of $B_1$
with an eigenvalue $\theta'$. Then the vector space $W'$ over
$\mbb{C}$ which is spanned by $(0,\mathbf{0}_2^t,\mathbf{z}^t)^t$ is
a thin irreducible $\mathcal{T}$-module over $\mbb{C}$ and dim
$W'=1$ if $\theta' \in \{ 0,-m\}$.
\end{enumerate}
 \end{lemma}
\begin{proof}
(1) Since $\mathbf{y}$ is the restricted eigenvector of $B_2$
associated with eigenvalue $\theta$, so $(B_2\mathbf{y})^t=(\theta
\mathbf{y})^t$. Hence, $(0,(B_2\mathbf{y})^t,\mathbf{0}_1^t)^t \in$
span$\{ (0,\mathbf{y}^t,\mathbf{0}_1^t)^t\}$, and $W$ is spanned by
$(0,\mathbf{y}^t,\mathbf{0}_1^t)^t$. Observe that
$A_1(0,\mathbf{y}^t,\mathbf{0}_1^t)^t=(0,(B_2\mathbf{y})^t,\mathbf{0}_1^t)^t$,
and also $\mathbf{y}$ is orthogonal to $\mathbf{1}_2$ with the
associated eigenvalue $0$ or $-m$. By Lemma \ref{res-eigen}
$L^t\mathbf{y}=\mathbf{0}_1$. Therefore, $W$ is $A_1$-invariant
and thus $\mc{M}$-invariant. $W$ is also $\mc{M}^*$-invariant.\\
(2) The proof is similar to that of (1).
\end{proof}

\begin{lemma}\label{decomp} For the association scheme $(K_{m})^{\wr 2}$,
let $V$ denote the standard $\mc{T}$-module. There
exist irreducible $\mc{T}$-modules $\{W_i\}_{1 \leq i \leq m(m-1)}$
and $\{W_j^{'}\}_{2 \leq j \leq m-2}$ such that $V = (\oplus
W_i)\oplus (\oplus W_j^{'})$.
\end{lemma}
\begin{proof} Using Gram Schmidt process we can find eigenvectors
$\mathbf{v}_1,\mathbf{v}_2,\cdots ,\mathbf{v}_{m(m-1)}$ of $B_2$
such that $\{\mathbf{v}_i\}_{1 \leq i \leq m(m-1)}$ span $E_1^{*}V
\cong \mbb{C}^{m(m-1)}$, $\langle \mathbf{v}_i,\mathbf{v}_j \rangle
=0$ for $i \neq j$ with $\mathbf{v}_1$ being the all-ones vector
$\mathbf{1}_2$. Let $\theta_i$ be the eigenvalue of $B_2$ with
respect to the eigenvector $\mathbf{v}_i$ for $2 \leq i \leq
m(m-1)$. Then $\theta_i \in \{0,-m\}$ for $2 \leq i \leq m(m-1)$.
Let $W_i$ denote the linear span of $(\langle
\mathbf{v}_i,\mathbf{v}_1 \rangle,\mathbf{0}_2^t,\mathbf{0}_1^t)^t$,
$(0,\mathbf{v}_i^{t},\mathbf{0}_1^t)^{t}$, and
$(0,\mathbf{0}_2^t,(L^t\mathbf{v}_i)^t)^{t}$ over $\mbb{C}$. Then
$W_i$ is a thin irreducible $\mathcal{T}$-module and $W_i \cap W_j =
\{\mathbf{0}\}$ for $i \neq j$. Also,
\begin{displaymath}
dim W_{i}= \left\{\begin{array}{ll}
3 & \textrm{if  } i = 1 \\
1 & \textrm{if  } 2 \leq i \leq m(m-1)\\
\end{array} \right.
\end{displaymath}
Note that $W_1$ is the primary module generated by $(\langle
\mathbf{v}_1,\mathbf{v}_1 \rangle,\mathbf{0}_2^t,\mathbf{0}_1^t)^t$,
$(0,\mathbf{v}_1^{t},\mathbf{0}_1^t)^{t}$, and
$(0,\mathbf{0}_2^t,(L^t\mathbf{v}_1)^t)^{t}$ over $\mbb{C}$. For
each $i$, $2\le i\le m(m-1)$, $W_i$ is generated by
$(0,\mathbf{v}_i^{t},\mathbf{0}_1^t)^{t}$.  Note that
$\mathbf{w}_1=L^t\mathbf{v}_1$ is an eigenvector of $B_1$. Let
$\mathbf{w}_2, \cdots, \mathbf{w}_{m-1}$ be the eigenvectors of
$B_1$ such that $\mathbf{w}_1, \cdots, \mathbf{w}_{m-1}$ span
$E_2^{*}V \cong \mathbb{C}^{m-1}$ and $\langle
\mathbf{w}_i,\mathbf{w}_j \rangle = 0$ for $ i\neq j$. Let $W_i^{'}$
be the linear span of $(0,\mathbf{0}_2^t,\mathbf{w}_i^{t})^t$ over
$\mbb{C}$ for $2 \leq i \leq m-1$. $W_i^{'}$ is a thin irreducible
$\mathcal{T}$-module of dimension $1$. Thus, we have $V = (\oplus
W_i)\oplus (\oplus W_j^{'})$ as desired.
\end{proof}

\begin{lemma}\label{irr-modules} For $(K_{m})^{\wr 2}$,
let $\{\theta_i\}_{1 \leq i \leq m(m-1)}$, $\{\theta_i^{'}\}_{2 \leq
i \leq m-1}$, $\{W_i\}_{1 \leq i \leq m(m-1)}$ and $\{W_i^{'}\}_{2
\leq i \leq m-1}$. Then the following hold.
\begin{enumerate}
\item For all $i$ with $2 \leq i \leq m(m-1)$, $W_1$ and $W_i$ are
not $\mathcal{T}$-isomorphic.
\item For all $i$ and $j$ with $1 \leq i \leq m(m-1)$ and
$2 \leq j \leq m-1$, $W_i$ and $W_j^{'}$ are not
$\mathcal{T}$-isomorphic.
\item For $i$ and $j$ with $2 \leq i,j \leq m(m-1)$, $W_i$ and $W_j$
are $\mathcal{T}$-isomorphic if and only if $\theta_i = \theta_j$.
\item For $i$ and $j$ with $2 \leq i,j \leq m-1$, $W_i^{'}$ and
$W_j^{'}$ are $\mathcal{T}$-isomorphic.
\end{enumerate}
\end{lemma}
\begin{proof} It is straightforward from the previous lemma
according to Lemma 3.4 in \cite{TY94}.
\end{proof}
\begin{remark} In order to describe the Terwilliger algebra of
$(K_m)^{\wr 2}$, let $\Lambda$ denote the index set for the
isomorphism classes of irreducible $\mc{T}$-modules. Following the
results of Lemmas \ref{decomp} and \ref{irr-modules}, we can group
the isomorphic irreducible modules together, say $V_{\lambda}$ for
$\lambda\in \Lambda$, so that
$V=\bigoplus\limits_{\lambda\in\Lambda}V_{\lambda}$. By Lemma
\ref{central}, we know that for each subspace $V_{\lambda}$ there is
a unique central idempotent $e_{\lambda}$ such that
$V_{\lambda}=e_{\lambda}V$. Let $W$ be an irreducible
$\mc{T}$-module in the decomposition of $V$. Then the map taking
$A\in e_{\lambda}\mc{T}$ to the endomorphism $\mathbf{w}\mapsto
A\mathbf{w}$ where $\mathbf{w}\in W$ is an isomorphism. Hence we
have $e_{\lambda}V\cong End_{\mbb{C}}W$. Thus
$\mc{T}=\bigoplus\limits_{\lambda\in \Lambda}e_{\lambda}\mc{T}$ is
isomorphic to a direct sum of complex matrix algebra $M_k(\mbb{C})$
where $k=dim(W)$ as in the following. In what follows we use the
notation $M_1(\mbb{C})^{\oplus l}$ for the direct sum $M_1(\mbb{C})
\oplus M_1(\mbb{C})\oplus \cdots \oplus M_1(\mbb{C})$ of $l$ copies
of $M_1(\mbb{C})$.
\end{remark}

\begin{theorem} Let $\mc{T}$ be the Terwilliger algebra of
$(K_{m})^{\wr 2}$. Then dim $\mathcal{T}= 12$ and $$\mathcal{T}
\cong M_3(\mbb{C}) \oplus M_1(\mbb{C})^{\oplus 3}.$$
\end{theorem}
\begin{proof} By Lemma \ref{irr-modules}, the list of non-isomorphic
irreducible $\mc{T}$-modules of $(K_m)^{\wr 2}$ consists of (i) the
primary module $W_1$ of dimension 3, (ii) two one-dimensional
non-isomorphic irreducible modules from $W_i$, $2\le i\le m(m-1)$,
that represent two isomorphism classes corresponding to the
eigenvalues $0$ and $-m$, and (iii) one one-dimensional irreducible
module which representing the class of $W_j^{\prime}$ for all $2\le
j\le m-1$.
\end{proof}

\subsection{The Terwilliger algebra of $(K_m)^{\wr 3}$ and its irreducible modules}

We now extend the irreducible modules of $(K_m)^{\wr 2}$ to find
irreducible modules of $(K_m)^{\wr 3}$ in this section, and then
generalize it to describe the Terwilliger algebra for $(K_m)^{\wr
d}$ for an arbitrary $d\ge 3$ in the next section. We note that the
Terwilliger algebra of $(K_m)^{\wr 2}$ is generated by $A_1$,
$E_0^{*}$ and $E_1^{*}$. As we move on to three or higher wreath
powers of $K_m$, the scheme is no longer $P$-polynomial, so, we must
consider many more, but within $2d+1$ generators for a $d$-power
case. We will see that a `concrete Terwilliger algebra' for the
wreath powers can be described as well. (Here the term `concrete' is
used in comparison with an `abstract' Terwilliger algebra described
in terms of generators and relations as in \cite{Eg00}.) We describe
the irreducible $\mc{T}$-modules of $(K_m)^{\wr d}$ by extending the
irreducible $\mc{T}$-modules of $(K_m)^{\wr (d-1)}$ for $d\ge 3$.
This can be done by investigating the structural relations between
the wreath square and the wreath cube to begin the iterative
process.

Let $\mathcal{X}= (K_{m})^{\wr 3}=(X, \{R_i\}_{0\le i\le 3})$ be a
$3$-class association scheme of order $m^3$. Let $X=\{(a, b, c): a,
b, c\in [m] \}$. Choose $x_1=(1,1,1)$ and fix it as the base vertex
and we will refer to it as $x$ henceforth. Let us order the
association relations such that $R_0(x)=\{(1,1,1) \}$, and
\[R_1(x) = \{(a, b, c):\ a, b\in [m],\ c\in \{2, 3, \dots, m\}\},\]
\[R_2(x)=\{ (a, b, 1):\ a\in [m],\ b\in \{2, 3, \dots, m\} \},\]
\[R_3(x)=\{ (a, 1, 1):\ a\in \{2, 3, \dots, m\} \}.\]
Let the relation table of $\mathcal{X}$ be decomposed into block
matrix form according to the partition $X= R_0(x) \cup R_1(x) \cup
R_2(x) \cup R_3(x)$. We can see that the relation matrix of
$(K_m)^{\wr 2}$ is embedded into that of $(K_m)^{\wr 3}$ by
relabeling the association relations. To see this, consider the
relation matrix $R$ for $(K_m)^{\wr 2}$ given in the previous
subsection. Denote the block labeled by $B_2+2C_2$ of $R$ by $D$.
Then $D$ has entries 0, 1, or 2. We now form a new block from $D$ by
replacing entries $2$ with $3$, $1$ with $2$ and $0$ with $0$, and
call this $D_2$. It is easy to see that $D_2$ is the block of the
relation table $(K_m)^{\wr 3}$ that is indexed by the vertices in
$R_2(x)$. The relation table of $(K_m)^{\wr 3}$ is given by

\[\left[\begin{array}{cccc}
0& \mathbf{1}_3^t & 2\mathbf{1}_2^t & 3\mathbf{1}_1^t \\
\mathbf{1}_3& D_3 & L_2 & L_1\\
2\mathbf{1}_2& L_2^{t} & D_2 & 2L\\
3\mathbf{1}_1 & L_1^{t} & 2L^{t} & 3C_1\\
\end{array}\right]
\]
where $\mathbf{1}_3,\ \mathbf{1}_2,\ \mathbf{1}_1$ are all-ones
column vectors of dimension $m^2(m-1),\ m(m-1),\ m-1$, respectively;
$L_2, L_1, L$ are all-ones matrices of size $m^{2}(m-1) \times
m(m-1),\ m^{2}(m-1) \times (m-1),\ m(m-1)\times (m-1)$,
respectively; $D_3 = I_m \otimes D_2 + (J_m-I_m) \otimes J_{m(m-1)}$
and $C_1= J_{m-1} - I_{m-1}$. We are now ready to see the
descriptions of the irreducible $\mc{T}$-modules of the wreath cube.

\begin{theorem} The primary irreducible $\mc{T}$-module of
$\mathcal{X}=(K_{m})^{\wr 3}$ is spanned by the four vectors
\[(1, \mathbf{0}_3^t,\mathbf{0}_2^t, \mathbf{0}_1^t)^t,\quad
(0, \mathbf{1}_3^t,\mathbf{0}_2^t, \mathbf{0}_1^t)^t,\quad (0,
\mathbf{0}_3^t,\mathbf{1}_2^t, \mathbf{0}_1^t)^t,\quad (0,
\mathbf{0}_3^t,\mathbf{0}_2^t, \mathbf{1}_1^t)^t \] where
$\mathbf{0}_3,\mathbf{0}_2, \mathbf{0}_1$ are all-zeros vectors of
dimensions $m^2(m-1),\ m(m-1),\ (m-1)$, respectively; and
$\mathbf{1}_3,\mathbf{1}_2, \mathbf{1}_1$ are all-ones vectors of
dimensions $m^2(m-1),\ m(m-1),\ (m-1)$, respectively.
\end{theorem}
\begin{proof} It is straightforward. \end{proof}
\begin{theorem} Let $V$ be the standard module of $(K_{m})^{\wr 3}$.
There exist irreducible $\mc{T}$-modules $\{W_i^{1}\}_{2 \leq i \leq
(m-1)}$ such that $(0, \mathbf{0}_3^t,\mathbf{0}_2^t,
\mathbf{1}_1^t)^t $ and  $\{W_i^{1}\}_{2\le i\le m-1}$ together
constitute $E_3^*V$.
\end{theorem}

\begin{proof} We have seen in the proof of Lemma \ref{decomp} that
there exist vectors $\mathbf{w}_1, \cdots, \mathbf{w}_{m-1}$  which
span $\mbb{C}^{m-1}$ and $\langle \mathbf{w}_i,\mathbf{w}_j \rangle
= 0$ for $1 \leq i\neq j \leq m-1$. Let $W_i^{1}$ be the linear span
of $(0,\mathbf{0}_3^{t},\mathbf{0}_2^{t},\mathbf{w}_i^{t})^t$ over
$\mbb{C}$. For $(K_m)^{\wr 2}$, the generators $E_i^{*}A_hE_i^{*}$
act on the modules $W_i^{'}$ and the significant nonzero actions are
$E_2^{*}A_0E_2^{*}$ and $E_2^{*}A_2E_2^{*}$. The embedded structure
of $(K_m)^{\wr 2}$ in $(K_m)^{\wr 3}$ ensures that for $(K_m)^{\wr
3}$ the generators $E_i^{*}A_hE_i^{*}$ act on the linear space
$W_i^{1}$, $2 \leq i \leq (m-1)$ and the only nonzero actions are
$E_3^{*}A_0E_3^{*}$ and $E_3^{*}A_3E_3^{*}$. It is clear that
$(0,\mathbf{0}_3^{t},\mathbf{0}_2^{t},\mathbf{w}_i^{t})^t$ are
$E_0^*A_3E_3^*, E_1^*A_1E_3^*, E_2^*A_2E_3^*$-invariant, and each
$W_i^{1}$ is an irreducible $\mathcal{T}$-module of dimension $1$.
Thus, the result follows.
\end{proof}
\begin{theorem} In the standard module $V$ of $(K_m)^{\wr 3}$,
there exist irreducible modules $\{W_i^{2}\}_{1 \leq i \leq m(m-1)}$
such that $(0, \mathbf{0}_3^t, \mathbf{1}_2^t, \mathbf{0}_1^t)^t$
and $\{W_i^2\}_{2\le i\le m(m-1)}$ constitute $E_2^{*}V$.
\end{theorem}

\begin{proof} We have seen in Lemma \ref{decomp} that there exist vectors
$\mathbf{v}_1, \dots, \mathbf{v}_{m(m-1)}$  such that $\mathbf{v}_1,
\dots, \mathbf{v}_{m(m-1)}$ span $\mbb{C}^{m(m-1)}$ and $\langle
\mathbf{v}_i,\mathbf{v}_j \rangle = 0$ for $1 \leq i\neq j \leq
m(m-1)$. Let $W_i^{2}$ be the linear span of
$(0,\mathbf{0}_3^{t},\mathbf{v}_i^{t},\mathbf{0}_1^{t})^t$ over
$\mbb{C}$. For $(K_m)^{\wr 2}$ the generators $E_i^{*}A_hE_i^{*}$
act on the modules $W_i$ and the significant nonzero actions are
$E_1^{*}A_0E_1^{*}$, $E_1^{*}A_1E_1^{*}$ and $E_1^{*}A_2E_1^{*}$.
The embedded structure of $(K_m)^{\wr 2}$ in $(K_m)^{\wr 3}$ ensures
that for $(K_m)^{\wr 3}$ the generators $E_i^{*}A_hE_i^{*}$ act on
the linear spaces $W_i^{2}$, $2 \leq i \leq m(m-1)$ and the nonzero
actions are  $E_1^{*}A_0E_1^{*}$, $E_2^{*}A_2E_2^{*}$ and
$E_2^{*}A_3E_2^{*}$. Also, we see that $(0, \mathbf{0}_3^t,
\mathbf{v}_1^t,\mathbf{0}_1^t)$ is invariant under the action of
$E_0^{*}A_2E_2^{*}$, $E_1^{*}A_1E_2^{*}$ and $E_3^{*}A_2E_2^{*}$.
Each of $\{W_i^{2}\}_{2 \leq i \leq m(m-1)}$ is a irreducible
$\mathcal{T}$-module of dimension $1$, and the result follows.
\end{proof}

We now describe the irreducible modules that span $E_1^{*}V$. We
know that $E_1^{*}V \cong \mbb{C}^{m^2(m-1)}$. We observe that any
$m^2(m-1)$ dimensional column vector can be partitioned into $m$
equal parts each with $m(m-1)$ components. Let each part be denoted
by the index $j$ where $j \in \{1,2, \dots ,m\}$ so that the
$m^2(m-1)$ dimensional column vector is of the form
$(\mathbf{u}_1^{t},\mathbf{u}_2^{t},\dots,\mathbf{u}_m^{t})^{t}$
where each $\mathbf{u}_j$ is an $m(m-1)$ dimensional column vector
for each $j \in \{1,2, \dots ,m\}$. Let
$\mathbf{u}_{j,i}=(\mathbf{u}_1^{t},\mathbf{u}_2^{t},\dots,
\mathbf{u}_l^{t},\dots,\mathbf{u}_m^{t})^{t}$ denote the $m^2(m-1)$
dimensional column vector such that $\mathbf{u}_l=\delta
_{lj}\mathbf{v}_i$ where $\mathbf{v}_i$ are the ones given in Lemma
\ref{decomp} and $\delta_{lj}$ is the Kronecker delta; i.e.,
$\delta_{lj}$ is $1$ if and only if $l= j$, and it is zero
otherwise.

\begin{lemma} In the standard module $V$ of $(K_{m})^{\wr 3}$,
for each pair $j, i$, $\ 1\le j\le m,\ 2 \leq i \leq m(m-1)$, the
linear subspace $W_{j,i}$ of spanned by
$(0,\mathbf{u}_{j,i}^{t},\mathbf{0}_2^{t},\mathbf{0}_1^{t})^t$ is an
irreducible $\mc{T}$-module contained in $E_1^*V$.
\end{lemma}

\begin{proof} In the Terwilliger algebra of $(K_m)^{\wr 2}$,
the generators $E_i^{*}A_hE_i^{*}$ act on the modules $W_i$ and the
only nonzero actions are due to $E_1^{*}A_0E_1^{*}$,
$E_1^{*}A_1E_1^{*}$ and $E_1^{*}A_2E_1^{*}$. For $(K_m)^{\wr 3}$,
the submatrix $D_3 = I_m \otimes D_2 + (J_m-I_m) \otimes J_{m(m-1)}$
indicates that in the Terwilliger algebra of $(K_m)^{\wr 3}$ the
generators $E_i^{*}A_hE_i^{*}$ act on the linear spaces in
$\{W_{j,i}:\ 1\le j\le m,\ 2 \leq i \leq m(m-1) \}$ where $W_{j,i}$
is the linear span of
$(0,\mathbf{u}_{j,i}^{t},\mathbf{0}_2^{t},\mathbf{0}_1^{t})^t$ and
the only nonzero actions are due to $E_1^{*}A_0E_1^{*}$,
$E_1^{*}A_2E_1^{*}$, $E_1^{*}A_3E_1^{*}$, $E_0^{*}A_1E_1^{*}$,
$E_2^{*}A_1E_1^{*}$ and $E_3^{*}A_1E_1^{*}$. It follows that
$W_{j,i}$ is an irreducible module of dimension $1$ for each $j,i$.
\end{proof}

\begin{remark}
We note that if $W^{0}$ denotes the linear span of
$(1,\mathbf{0}_3^{t},\mathbf{0}_2^{t},\mathbf{0}_1^{t})^t$ over
$\mbb{C}$, then $W^{0}$ is an irreducible module that spans
$E_0^{*}V$. We recall that, by Lemma \ref{irr-modules}, among the
irreducible modules $\{W_{i}\}_{2 \leq i \leq m(m-1)}$ of
$(K_m)^{\wr 2}$ there are two non-isomorphic $\mathcal{T}$-modules
of dimension $1$. In sum, so far, for $(K_m)^{\wr 3}$ we have the
following non-isomorphic $\mathcal{T}$-modules
\begin{enumerate}
\item The primary module of dimension $4$.
\item Two non-isomorphic $\mathcal{T}$-modules of dimension $1$ in $E_1^{*}V$.
\item Two non-isomorphic $\mathcal{T}$-modules of dimension $1$ in $E_2^{*}V$.
\item One non-isomorphic $\mathcal{T}$-modules of dimension $1$ in $E_3^{*}V$.
\end{enumerate}
\end{remark}

It leaves us with one non-isomorphic irreducible
$\mathcal{T}$-module of dimension $1$ which is not accounted for as
the total dimension of the Terwilliger algebra must be 22 by the
Theorem \ref{dimension}. It is in the subconstituent $E_1^{*}V$ as
in the following.

\begin{lemma} Pick any nonzero vector $\mathbf{u} \in E_1^{*}V$ which is
orthogonal to all the irreducible modules $\{W_{j,i}:\ 1\le j\le m,\
2 \leq i \leq m(m-1)\}$. The $\mathbf{u}$ spans a one-dimensional
irreducible $\mathcal{T}$-module. The actions of
$E_1^{*}A_0E_1^{*}$, $E_1^{*}A_1E_1^{*}$, $E_1^{*}A_2E_1^{*}$,
$E_1^{*}A_3E_1^{*}$, $E_0^{*}A_1E_1^{*}$, $E_2^{*}A_1E_1^{*}$ and
$E_3^{*}A_1E_1^{*}$ on this $\mc{T}$-module are nonzero and rest are
all zero.
\end{lemma}
\begin{proof} It is straightforward.
\end{proof}

\begin{theorem} For $(K_{m})^{\wr 3}$,
 $\mathcal{T} \cong M_4(\mbb{C}) \oplus M_1(\mbb{C})^{\oplus 6}$.
\end{theorem}

\begin{proof} It follows from the above remark
and lemmas.
\end{proof}

\subsection{The Terwilliger algebra of $(K_m)^{\wr d}$ for $m \geq
3$}

The description of the concrete Terwilliger algebra of a scheme
involves describing the irreducible modules that constitute
different subconstituents of the algebra. Earlier in this section,
we have already studied all the subconstituents of the $2$-class and
$3$-class wreath power association schemes. From the irreducible
modules of $(K_m)^{\wr 3}$ we can describe all irreducible modules
for $(K_m)^{\wr 4}$, from $(K_m)^{\wr 4}$ we can describe all
irreducible modules for $(K_m)^{\wr 5}$, and so on. We now develop a
recursive method to describe the Terwilliger algebra of a $d$-class
association scheme $(K_m)^{\wr d}$ from a $(d-1)$-class scheme
$(K_m)^{\wr (d-1)}$. Let $\mathcal{X}=(X, \{R_i\}_{0 \leq i \leq
d})$ denote the $d$-class scheme $(K_m)^{\wr d}$ with \[X = [m]
\times [m] \times \cdots \times [m]= \{(a_1,a_2,\cdots,a_d):\ a_i
\in [m], \mbox{ for } i = 1,2,\dots,d\}\] Let $(1,1,\dots,1) \in X$
be a fixed base vertex $x$ of $\mathcal{X}$. Without loss of
generality, we can arrange the label of associate relations and the
vertices so that for $i=1,2,\dots,d-1$,
    \[R_i(x)=\{(a_1,a_2,\dots,a_{d-i-1},a,1,1,\dots,1):\ a_k
    \in [m] \mbox{ for } 1\le k\le d-i-1,\ a \in
    [m]-\{1\}\};\]
    \[R_0(x)=\{x\};\]
\[R_d(x)=\{(a,1,1,\dots,1):~a \in [m]-\{1\}\}.\]
In the same manner, we can get the relation table for $d-1$ wreath
power as well. First let us look at the relation table of
$(K_m)^{\wr (d-1)}$.

\begin{displaymath}
\left[\begin{array}{cccccc}
0             &\mathbf{1}_{d-1}^t      &2\mathbf{1}_{d-2}^t     & \cdots
&(d-2)\mathbf{1}_{2}^t  &(d-1)\mathbf{1}_{1}^t  \\
\mathbf{1}_{d-1}   &T_{d-1}      & J_{d-1,d-2}      &\cdots     &J_{d-1,2}    & J_{d-1,1}  \\
2\mathbf{1}_{d-2}  & J_{d-1,d-2}^{t}  & T_{d-2}     & \cdots    & 2J_{d-2,2}  & 2J_{d-2,1}\\
      \vdots    &   \vdots &  \vdots      &  \ddots     &  \vdots     & \vdots \\
(d-2)\mathbf{1}_{2}& J_{d-1,2}^{t}& 2J_{d-2,2}^{t}&  \cdots   & T_{2}      & (d-2)J_{1,1}\\
(d-1)\mathbf{1}_{1}& J_{d-1,1}^{t} & 2J_{d-2,1}^{t}&  \cdots  & (d-2)J_{1,1}^{t}  &T_{1} \\
\end{array}\right]
\end{displaymath}
where $\mathbf{1}_i$ are all-ones column vectors of size
$m^{i-1}(m-1)$, $J_{j,k}$ are all-ones matrices of size
$m^{j-1}(m-1)\times m^{k-1}(m-1)$, $T_1=(d-1)(J_{m-1}-I_{m-1})$ and
$T_i=I_m \otimes T_{i-1}+(d-i)(J_m-I_m)\otimes J_{m^{i-2}(m-1)}$ for
$i \in \{2, 3,\dots,d-1\}$.

Now using this table, we can describe the relation table for
$(K_m)^{\wr d}$ as follows. In the diagonal blocks $T_i$ make the
following changes. The entry $0$ is kept same, and the entries $i$
are replaced with $i+1$, for all $i=1, 2, \dots, d-1$. Let us name
the resulting new blocks $U_i$ for $i \in \{1,2,\dots,d-1\}$. It is
not hard to see that $U_i$ are the diagonal blocks of the relation
table of $(K_m)^{\wr d}$. Let $U_d=I_m \otimes
U_{d-1}+(J_m-I_m)\otimes J_{m^{d-2}(m-1)}$, $\mathbf{1}_d$ is the
all-ones column vector of size $m^{d-1}(m-1)$. Abusing notation and
denoting all the all-ones matrices in the relation table as $J$ for
all dimensions the relation table of the $d$-class association
scheme is
 \begin{displaymath}
\left[\begin{array}{ccccccc} 0         & \mathbf{1}_{d}^t
&2\mathbf{1}_{d-1}^t      &3\mathbf{1}_{d-2}^t     & \cdots
&(d-1)\mathbf{1}_{2}^t  & d\mathbf{1}_{1}^t  \\
\mathbf{1}_{d}&  U_{d}  & J_{d,d-1}   &   J_{d, d-2}   &  \cdots   & J_{d, 2}   & J_{d,1} \\
2\mathbf{1}_{d-1} & J_{d,d-1}^t &U_{d-1}      & 2J_{d-1,d-2}      &\cdots     &2J_{d-1,2}
 & 2J_{d-1,1} \\
3\mathbf{1}_{d-2} &J_{d,d-2}^t & 2J_{d-1,d-2}^t  & U_{d-2}     & \cdots    &3 J_{d-2,2}
& 3J_{d-2,1}\\
      \vdots  &\vdots  &   \vdots &  \vdots      &  \ddots     &  \vdots     & \vdots \\
(d-1)\mathbf{1}_{2}& J_{d,2}^t &2J_{d-1,2}^t & 3J_{d-2,2}^t&  \cdots   & U_{2}      & (d-1)J_{1,1}\\
d\mathbf{1}_{1}&J_{d,1}^t &2J_{d-1,1}^t & 3J_{d-2,1}^t&  \cdots  & (d-1)J_{1,1}^t  &U_{1} \\
\end{array}\right]
\end{displaymath}

In the next couple of paragraphs we will discuss the subconstituents
of the $d$-class association scheme $(K_m)^{\wr d}$. Let
$\mathbf{0}_{i}$ denote zero column vectors of size $m^{i-1}(m-1)$
for $1 \leq i \leq d$. Any $m^d$ dimensional column vectors can be
divided into subparts $1$, $m^{d-1}(m-1)$, $m^{d-2}(m-1)$, $\dots$,
$m(m-1)$ and $(m-1)$ respectively so that any vector looks like
$(p,\mathbf{p}_d^{t},\cdots,\mathbf{p}_1^{t})^{t}$ where
$\mathbf{p}_i$ is a $m^{i-1}(m-1)$ dimensional column vector for $1
\leq i \leq d$. With these notations, the primary module of the
$d$-class association scheme $(K_m)^{\wr d}$ may be described as
follows.

\begin{theorem} Let $V$ be the standard module of $(K_m)^{\wr d}$. The
vector
 $(1,\mathbf{0}_{d}^{t},\dots,\mathbf{0}_{1}^{t})^{t}$ and vectors
$\mathbf{q}_i=(0,\mathbf{p}_d^{t},\dots,\mathbf{p}_j^{t},\dots,\mathbf{p}_1^{t})^{t}$
for $1 \leq i \leq d$ such that
\begin{displaymath}
\mathbf{p}_j= \left\{\begin{array}{ll}
\mathbf{1}_j & \textrm{if  } i=j \\
\mathbf{0}_{j} & \textrm{if  }i \neq j \\
\end{array} \right.
\end{displaymath}
generates the primary $\mathcal{T}$-module.
\end{theorem}
\begin{proof} Straightforward.
\end{proof}

Let us consider the $d$-class association scheme $(K_{m})^{\wr d}$
of order $m^d$ and let $V$ denote its standard module. Finding the
irreducible modules of the subconstituents $E_{i}^{*}V$ for $2 \leq
i \leq d$ and $E_{0}^{*}V$ is more routine. $E_1^{*}V$ needs to be
treated differently than the other and we will come to that as we go
along.

Suppose that the $m^{d-1}$ dimensional column vectors
$(0,\mathbf{h}_j^{t})^{t}$ for $1 \leq j \leq m^{d-i}(m-1)-1$
generate the one dimensional modules of the subconstituent
$E_{i-1}^{*}V$ for the $(d-1)$-class association scheme $(K_m)^{\wr
(d-1)}$. If we add the $m^{d-1}(m-1)$ dimensional column vector
$\mathbf{0}_{d}^{t}$ right after $0$ in the above vectors we land up
with $m^d$ dimensional vectors. For $1 \leq j \leq m^{d-i}(m-1)-1$
the vectors $(0,\mathbf{0}_d^t, \mathbf{h}_j^{t})^{t}$ and
$(0,\mathbf{0}_{d}^{t},\dots,\mathbf{1}_{d-i-1}^{t},\dots,\mathbf{0}_{1}^{t})^{t}$
span $\mathbb{C}^{m^{d-i}(m-1)}$.  Also,
 $\langle (0,\mathbf{h}_j^{t}),(0,\mathbf{h}_k^{t}) \rangle = 0$ for
$j,k \in \{1,2,\dots, m^{d-i}(m-1)-1\}$.  For the scheme $(K_m)^{\wr
(d-1)}$ since $(0,\mathbf{h}_j^{t})^{t}$ generates a one dimensional
module it is $E_i^{*}A_jE_h^{*}$ invariant for all $i,j,h \in
\{0,1,\dots,d-1\}$. For $1 \leq j \leq m^{d-i}(m-1)-1$, let
$W_{j}^{\wr d}$ be the linear span of the vector
$(0,\mathbf{0}_{d}^{t},\mathbf{h}_j^{t})^{t}$. The embedded
structure of the $(d-1)$-class association scheme $(K_m)^{\wr
(d-1)}$ in the $d$-class scheme $(K_m)^{\wr d}$ ensures that for $1
\leq j \leq m^{d-i}(m-1)-1$,
$(0,\mathbf{0}_{d}^{t},\mathbf{h}_j^{t})^{t}$ are
$E_i^{*}A_jE_h^{*}$ invariant for all $i,j,h \in \{0,1,\dots,d\}$.
Note that now we are talking about the triple products of the
$d$-class scheme $(K_m)^{\wr d}$. The vectors
$(0,\mathbf{0}_{d}^{t},\mathbf{h}_j^{t})^{t}$ for $1 \leq j \leq
m^{d-i}(m-1)-1$
 generate all the one dimensional non-primary irreducible
modules of the subconstituent $E_{i}^{*}V$ for the $D$-class
association scheme.

We have so far described the subconstituents $E_0^{*}V$, $E_2^{*}V$,
$\dots$, $E_d^{*}V$ for the scheme $(K_m)^{\wr d}$. Now, $E_1^{*}V
\cong \mathbb{C}^{m^{d-1}(m-1)}$. Observe that any $m^{d-1}(m-1)$
dimensional column vector can be partitioned into $m$ equal parts
each of dimension $m^{d-2}(m-1)$. Let each part be denoted by the
index $j$ where $j \in \{1,2, \dots, m\}$ so that the $m^{d-1}(m-1)$
dimensional vector is of the form
$(\mathbf{r}_1^{t},\mathbf{r}_2^{t},\dots,\mathbf{r}_m^{t})^{t}$
where each $\mathbf{r}_j$, $j \in \{1,2, \dots, m\}$ is a
$m^{d-2}(m-1)$ dimensional column vector. For each $i \in
\{2,3,\dots, m^{d-2}(m-1)\}$ and $j \in \{1,2, \dots , m\}$ let
$\mathbf{r}_{j,i}=(\mathbf{r}_1^{t},\mathbf{r}_2^{t},\dots,
\mathbf{r}_l^{t},\dots,\mathbf{r}_m^{t})^{t}$ denote the
$m^{d-1}(m-1)$ dimensional column vector such that
$\mathbf{r}_l=\delta _{lj}\mathbf{h}_j$ where $(0,\mathbf{h}_j^{t})$
generates the modules of $E_1^{*}V$ for the $(d-1)$-class scheme
$(K_m)^{\wr (d-1)}$. It is easy to see that each of the vectors
$(0,\mathbf{r}_{j,i}^{t},\mathbf{0}_{d-1}^{t},\dots,\mathbf{0}_1^{t})^{t}$
generates one dimensional irreducible module in the subconstituent
$E_{1}^{*}V$ for the scheme $(K_m)^{\wr d}$. These do not constitute
all the irreducible modules in $E_1^{*}V$. So far we have accounted
for the following non-isomorphic $\mathcal{T}$-modules for
$(K_m)^{\wr d}$.
\begin{enumerate}
\item The primary module of dimension $d+1$.
\item $d-1$ non-isomorphic $\mathcal{T}$-modules of
    dimension $1$ in $E_1^{*}V$.
\item $d-1$ non-isomorphic $\mathcal{T}$-modules of
    dimension $1$ in $E_2^{*}V$.
\item $d-2$ non-isomorphic $\mathcal{T}$-modules of
    dimension $1$ in $E_3^{*}V$

         $\quad \vdots$

\item $2$ non-isomorphic $\mathcal{T}$-modules of dimension
    $1$ in $E_{d-1}^{*}V$.
\item $1$ non-isomorphic $\mathcal{T}$-modules of dimension
    $1$ in $E_d^{*}V$.
\end{enumerate}

From Theorem \ref{dimension} we have the dimension of
$\mathcal{T}((K_m)^{\wr d}$ is $(d+1)^2 + \frac12d(d+1)$.
 From our above discussion we have $(d+1)^{2}+(d-1)+(d-1)+(d-2)+\cdots+2+1$
  of the dimension. That leaves us with one non-isomorphic $\mathcal{T}$
 module of dimension $1$ in the subconstituent $E_1^{*}V$ of multiplicity $m-1$.
 Pick any nonzero vector $\mathbf{r} \in E_1^{*}V$
 which is orthogonal to all the irreducible modules.
 Then $\mathbf{r}$ spans a one-dimensional irreducible $\mathcal{T}$-module.

In the discussion above we saw how we could build the irreducible
modules of the scheme $(K_m)^{\wr d}$.  We will conclude the section
by collecting all our non-isomorphic $\mathcal{T}$-modules and
describing the Terwilliger algebra of $(K_m)^{\wr d}$. (Further
detailed explanation of the content of this section can be found in
\cite{Bh08}.)

\begin{theorem} Let $\mathcal{X}= (K_m)^{\wr d}$ be a $d$-class association
scheme of order $m^d$, $m\ge 3$. Then the dimension of $\mathcal{T}$
is $(d+1)^2 + \frac12d(d+1)$ and $\mathcal{T} \cong
M_{d+1}(\mathbb{C}) \oplus M_1(\mbb{C})^{\oplus\frac{1}{2}d(d+1)}.$
\end{theorem}

\subsection{The Irreducible $\mathcal{T}$-modules of $(K_2)^{\wr d}$}

We note that for the case $(K_2)^{\wr d}$ the number of nonzero
triple products $E_i^{*}A_jE_h^{*}$ are fewer in number than in the
general case $(K_m)^{\wr d}$ with $m \geq 3$. What is nice for
$(K_2)^{\wr d}$ is that instead of just knowing the existence of
vectors that generate the irreducible modules, we are actually able
to get specific vectors that generate the irreducible
$\mc{T}$-modules.

Let $\mathcal{X}=(K_2)^{\wr d}$ be a $d$-class association scheme of
order $2^d$. Without loss of generality we label the $2^d$ elements
of $X$ by $x_1,x_2,\dots, x_{2^d}$. We fix $x_1$ as the base vertex
and we will refer to it as $x$ henceforth. Then, without loss of
generality, we can arrange the elements such that $X$ is partitioned
with $ R_0(x) = \{x \}$, $R_1(x)=\{x_{2}, \dots, x_{2^{d-1}+1} \}$,
consisting of the $2^{d-1}$ elements, $R_2(x)=\{x_{2^{d-1}+2}, \dots
, x_{2^{d-1}+2^{d-2}+1}\}$ consisting of the next $2^{d-2}$
elements, and so on, with the last part $R_d(x)=\{ x_{2^d} \}$.

%\begin{example}The relation table for a $4$-class scheme $(K_2)^{\wr 4}$ is the following
%
%\[\left(\begin{array}{c|cccccccc|cccc|cc|c}
%
% 0 & 1 & 1& 1& 1& 1& 1& 1& 1& 2& 2& 2& 2& 3& 3& 4\\\hline
%1 & 0 &4 &3 & 3& 2& 2&2 &2 & 1& 1& 1&1 & 1& 1&1 \\
%1 & 4 & 0& 3&3 &2 & 2& 2& 2&1 &1 &1 &1 & 1&1 &1 \\
% 1&  3&3 &0 &4 &2 & 2& 2&2 &1 &1 & 1& 1& 1& 1& 1\\
% 1&3  &3 &4 &0 & 2& 2& 2&2 &1 & 1&1 & 1&1 & 1&1 \\
% 1& 2 & 2&2 & 2&0 &4 & 3& 3&1 & 1&1 & 1& 1& 1&1 \\
%  1& 2 & 2&2 &2 &4 &0 & 3&3 &1 &1 &1 & 1& 1&1 &1 \\
%  1 & 2 & 2&2 &2 & 3& 3& 0& 4& 1& 1& 1& 1&1 & 1&1 \\
%  1   & 2 & 2&2 & 2&3 &3 & 4&0 &1 & 1&1 & 1& 1&1 &1 \\\hline
%  2    & 1 & 1& 1& 1& 1&1 &1 &1 &0 &4 &3 & 3&2 &2 &2 \\
%   2     & 1 &1 & 1&1 & 1&1 &1 & 1&4 &0 & 3& 3&2 &2 & 2\\
%   2      & 1 & 1&1 & 1&1 & 1& 1&1 & 3& 3& 0& 4&2 &2 & 2\\
%    2      & 1 & 1&1 & 1&1 & 1&1 &1 &3 &3 & 4& 0&2 &2 &2 \\\hline
%    3      & 1 &1 &1 & 1& 1&1 &1 & 1&2 & 2& 2&2 &0 & 4&3 \\
%    3     & 1 & 1& 1&1 &1 &1 &1 &1 & 2&2 &2 & 2&4 & 0&3   \\\hline
%      4  & 1 & 1& 1& 1&1 &1 &1 &1 & 2& 2&2 &2 & 3&3 &0 \\
%\end{array}\right)\]
%\end{example}
%
%We will describe vectors that will generate the irreducible
%$\mathcal{T}$-modules of $(K_2)^{\wr d}$.

\begin{remark} Let $\mathbf{1}$ denote the all-ones vector in the standard
module. Then the vector space over $\mathbb{C}$ spanned by
$\{E_i^{*}\mathbf{1}:\ 0 \leq i \leq d\}$ is a thin irreducible
$\mathcal{T}$-module of dimension $d+1$ and is the primary module
denoted as $\mathcal{P}$; so, by setting
$\mathbf{d}_i=E_i^*\mathbf{1}$, the set $\{\mathbf{d}_i:\ 0 \leq i
\leq d\}$ generates $\mc{P}$.\end{remark}

%\subsection{Construction of some new vectors}\label{vec}

In order to describe the irreducible $\mathcal{T}$-modules, whose
orthogonal direct sum forms the standard module $V$, we employ a
particular set of vectors. Let $\hat{\mathbf{x}}$ denote the column
vector with $1$ in the $x$-th position and $0$ elsewhere.

\begin{lemma} \label{vect} For $l \in\{1,2, \dots , d-1\}$ define set of vectors
$\{\mathbf{d}_{i}^{l}\}_{1 \leq i \leq 2^{d-l}-1}$ by
\[\mathbf{d}_{i}^{l}=
\sum_{k=0}^{2^{l-1}-1}\hat{\mathbf{x}}_{i+j+k}-
\sum_{k=2^{l-1}}^{2^l-1}\hat{\mathbf{x}}_{i+j+k}.\] For each $i$,
the corresponding values of $j$ are successively
\[j=1,\ 1+2^l-1,\ 1+2(2^l-1),\ 1+3(2^l-1),\ \dots,\ 1+(2^{d-l}-2)(2^l-1).\]
Then
\begin{enumerate}
\item[]  $\{\mathbf{d}_{i}^{1}\}_{1 \leq i \leq
    2^{d-1}-1}$;\quad
    $\mathbf{d}_{i}^{1}=\hat{\mathbf{x}}_{2i}-\hat{\mathbf{x}}_{2i+1}$

\item[] $\{\mathbf{d}_{i}^{2}\}_{1 \leq i \leq 2^{d-2}-1}$;\quad
    $\mathbf{d}_{i}^{2}=\hat{\mathbf{x}}_{i+j}+\hat{\mathbf{x}}_{i+j+1}-
    \hat{\mathbf{x}}_{i+j+2}-\hat{\mathbf{x}}_{i+j+3}$;
    $j=1,4,7,\dots$

\item[] $\{\mathbf{d}_{i}^{3}\}_{1 \leq i \leq
    2^{d-3}-1}$;\quad $\mathbf{d}_{i}^{3}=
    \sum_{k=0}^{2^{3-1}-1}\hat{\mathbf{x}}_{i+j+k}-
    \sum_{k=2^{3-1}}^{2^3-1}\hat{\mathbf{x}}_{i+j+k}$;
    $j=1,8,15,\dots$
\item[] \quad $\vdots$
\item[]
    $\mathbf{d}_{1}^{d-1}=\sum_{i=2}^{2^{d-2}+1}\hat{\mathbf{x}}_{i}-
    \sum_{i=2^{d-2}+2}^{2^{d-1}+1}\hat{\mathbf{x}}_{i}$.
\end{enumerate}
%\begin{example} The vectors for $(K_2)^{\wr 4}$
%are \[\begin{array}{ll} (0,1,-1,0,0,0,0,0,0,0,0,0,0,0,0,0)^{t},
%&
%(0,0,0,1,-1,0,0,0,0,0,0,0,0,0,0,0)^{t},\\
%(0,0,0,0,0,1,-1,0,0,0,0,0,0,0,0,0)^{t},&
%(0,0,0,0,0,0,0,1,-1,0,0,0,0,0,0,0)^{t},\\
%(0,0,0,0,0,0,0,0,0,1,-1,0,0,0,0,0)^{t}, &
%(0,0,0,0,0,0,0,0,0,0,0,1,-1,0,0,0)^{t},\\
%(0,0,0,0,0,0,0,0,0,0,0,0,0,1,-1,0)^{t},&
%(0,1,1,-1,-1,0,0,0,0,0,0,0,0,0,0,0)^{t},\\
%(0,0,0,0,0,1,1,-1,-1,0,0,0,0,0,0,0)^{t},&
%(0,0,0,0,0,0,0,0,0,1,1,-1,-1,0,0,0)^{t},\\
%(0,1,1,1,1,-1,-1,-1,-1,0,0,0,0,0,0,0)^{t}. & \end{array} \]
%\end{example}
%
%\subsection{Decomposition of the standard module into irreducible $\mathcal{T}$-modules}
In particular, $\langle \mathbf{d}_{i}^{l},\mathbf{d}_{h}^{k}
\rangle = 0$ unless $h=i$ and $k=l$.
\end{lemma}
\begin{proof} Proof follows from the construction of the vectors.
\end{proof}
\begin{lemma} Let $W_{\mathbf{d}_i^l}$ be the linear span of $\mathbf{d}_{i}^{l}$ for all the
vectors described in Lemma \ref{vect}. Then $W_{\mathbf{d}_i^l}$ is
an irreducible $\mathcal{T}$-module of dimension $1$ for $(K_2)^{\wr
d}$.
\end{lemma}
\begin{proof} %To show that a nonzero subspace $W \subseteq \mathbb{C}^{|X|}$
% is an irreducible $\mathcal{T}$-module we need to show that
 % $BW \subseteq W$ for all $B \in  \mathcal{T}$ and $W$
  % contains no $\mathcal{T}$-modules other than $0$ and $W$.
   To prove that $W_{\mathbf{d}_i^l}$ is an irreducible $\mathcal{T}$-module,
    we look at the action of the nonzero generators
    $E_i^{*}A_jE_h^{*}$ of the Terwilliger algebra on the
    vector $\mathbf{d}_{i}^{l}$. We will consider three cases.
\begin{enumerate}
\item The case when $l=1$:
\begin{enumerate}
\item For $\{\mathbf{d}_{i}^{1}\}_{1 \leq i \leq
    2^{d-2}}$ by the construction of
    $\mathbf{d}_{i}^{1}$, the action of the triple
    products $E_i^{*}A_jE_h^{*}$ on
    $\mathbf{d}_{i}^{1}$ can be nonzero only for
    $E_1^{*}A_0E_1^{*},E_1^{*}A_1E_1^{*},\dots,E_1^{*}A_dE_1^{*}$
    and $E_i^{*}A_1E_1^{*}$ for $2 \leq i \leq d$.
    Among them $E_1^{*}A_1E_1^{*}=0$. The generators
    $E_i^{*}A_jE_h^{*}$ of $\mathcal{T}$ that act on
    $\mathbf{d}_{i}^{1}$ in a nonzero manner are
    $E_1^{*}A_0E_1^{*}$ and $E_1^{*}A_dE_1^{*}$ with
    \[E_1^{*}A_0E_1^{*}\mathbf{d}_{i}^{1} =
    \mathbf{d}_{i}^{1}, \quad
    E_1^{*}A_dE_1^{*}\mathbf{d}_{i}^{1} =
    -\mathbf{d}_{i}^{1}.\]

\item  For $\{\mathbf{d}_{i}^{1}\}_{2^{d-2}+1 \leq i
    \leq 2^{d-2}+ 2^{d-3}}$ by the construction of
    $\mathbf{d}_{i}^{1}$, the action of the triple
    products $E_i^{*}A_hE_h^{*}$ on
    $\mathbf{d}_{i}^{1}$ are nonzero only for
    $E_2^{*}A_0E_2^{*}$ and $E_2^{*}A_dE_2^{*}$, and
    act in the following manner:
\[E_2^{*}A_0E_2^{*}\mathbf{d}_{i}^{1} = \mathbf{d}_{i}^{1},\quad
E_2^{*}A_dE_1^{*}\mathbf{d}_{i}^{1} = -\mathbf{d}_{i}^{1}.\]

\item For every set of vectors corresponding to the
    different subconstituents will follow the same
    pattern. The last case is the following.
For $\mathbf{d}_{2^{d-1}-1}^{1}$ the generators $E_i^{*}A_jE_h^{*}$
of $\mathcal{T}$ that act on $\mathbf{d}_{2^{d-1}-1}^{1}$ in a
nonzero manner for
\[E_{d-1}^{*}A_0E_{d-1}^{*}\mathbf{d}_{i}^{1} =
\mathbf{d}_{i}^{1},\quad E_{d-1}^{*}A_dE_{d-1}^{*}\mathbf{d}_{i}^{1}
= \mathbf{d}_{i}^{1}.\]
\end{enumerate}

\item The case when $l \in \{2,3, \dots, d-2\}$:
\begin{enumerate}
\item For $\{\mathbf{d}_{i}^{l}\}_{1 \leq i \leq
    2^{d-(l+1)}}$ by the construction of
    $\mathbf{d}_{i}^{l}$, the action of the triple
    products $E_i^{*}A_jE_h^{*}$ on
    $\mathbf{d}_{i}^{1}$ can be nonzero only for
    $E_1^{*}A_0E_1^{*},E_1^{*}A_1E_1^{*},\dots,E_1^{*}A_dE_1^{*}$.
    Among them $E_1^{*}A_1E_1^{*}=0$. The generators
    $E_i^{*}A_jE_h^{*}$ of $\mathcal{T}$ that act on
    $\mathbf{d}_{i}^{l}$ in a nonzero manner are
    $E_1^{*}A_0E_1^{*}, E_1^{*}A_dE_1^{*}$,
    $E_1^{*}A_{d-1}E_1^{*}, \dots,
    E_1^{*}A_{d-(l-1)}E_1^{*}$ in the following manner:
\[E_1^{*}A_0E_1^{*}\mathbf{d}_{i}^{l} =
\mathbf{d}_{i}^{l},\quad E_1^{*}A_dE_1^{*}\mathbf{d}_{i}^{l} =
\mathbf{d}_{i}^{l},\quad E_1^{*}A_{d-j}E_1^{*}\mathbf{d}_{i}^{l} =
2^{j}\mathbf{d}_{i}^{l} \mbox{ where } 1 \leq j \leq l-1.\]

\item For $\{\mathbf{d}_{i}^{l}\}_{2^{d-(l+1)}+1
    \leq i \leq 2^{d-(l+1)}+ 2^{d-(l+2)}}$ by the
    construction of $\mathbf{d}_{i}^{1}$, the action
    of the triple products $E_i^{*}A_jE_h^{*}$ on
    $\mathbf{d}_{i}^{1}$ can be nonzero only for
    $E_2^{*}A_0E_2^{*},E_2^{*}A_1E_2^{*},\dots,E_2^{*}A_dE_2^{*}$.
    Among them $E_2^{*}A_1E_2^{*}=0$. The generators
    $E_i^{*}A_jE_h^{*}$ of $\mathcal{T}$ that act on
    $\mathbf{d}_{i}^{l}$ in a nonzero manner are
    $E_2^{*}A_0E_2^{*}$, $E_2^{*}A_dE_2^{*}$,
    $E_2^{*}A_{d-1}E_2^{*}, \dots,
    E_2^{*}A_{d-(l-1)}E_2^{*}$ in the following manner:
\[E_2^{*}A_0E_2^{*}\mathbf{d}_{i}^{l} =
\mathbf{d}_{i}^{l},\quad E_2^{*}A_dE_2^{*}\mathbf{d}_{i}^{l} =
\mathbf{d}_{i}^{l},\quad E_2^{*}A_{d-j}E_2^{*}\mathbf{d}_{i}^{l} =
2^{j}\mathbf{d}_{i}^{l} \mbox{ for } 1 \leq j \leq l-2,\] and
\[E_2^{*}A_{d-(l-1)}E_2^{*}\mathbf{d}_{i}^{l} =-
2^{l-1}\mathbf{d}_{i}^{l}.\]

\item Following a similar reasoning the last case will
    be as follows. For $\mathbf{d}_{2^{d-l}-1}^{l}$ the generators $E_i^{*}A_jE_h^{*}$
of $\mathcal{T}$ that act on $\mathbf{d}_{2^{d-l}-1}^{l}$ in a
nonzero manner are $E_{d-1}^{*}A_0E_{d-1}^{*}$ and
$E_{d-1}^{*}A_dE_{d-1}^{*}$ in the following manner:
\[E_{d-l}^{*}A_0E_{d-l}^{*}\mathbf{d}_{2^{d-l}-1}^{l} =
\mathbf{d}_{2^{d-l}-1}^{l},\quad
E_{d-l}^{*}A_dE_{d-l}^{*}\mathbf{d}_{2^{d-l}-1}^{l} =
\mathbf{d}_{2^{d-l}-1}^{l},\]
\[E_{d-l}^{*}A_{d-j}E_{d-l}^{*}\mathbf{d}_{2^{d-l}-1}^{l} =
2^{j}\mathbf{d}_{2^{d-l}-1}^{l} \mbox{ for } 1 \leq j \leq l-2,\]
and
\[E_{d-l}^{*}A_{d-(l-1)}E_{d-l}^{*}\mathbf{d}_{2^{d-l}-1}^{l} =-
2^{l-1}\mathbf{d}_{2^{d-l}-1}^{l}.\]
\end{enumerate}
\item The case when $l=d-1$:

By the construction of $\mathbf{d}_{1}^{d-1}$, the action of the
triple products $E_i^{*}A_jE_h^{*}$ on $\mathbf{d}_{1}^{d-1}$ can be
nonzero only for
$E_1^{*}A_0E_1^{*},E_1^{*}A_1E_1^{*},\dots,E_1^{*}A_dE_1^{*}$. Among
them $E_1^{*}A_1E_1^{*}=0$. The generators $E_i^{*}A_jE_h^{*}$ of
$\mathcal{T}$ that act on $\mathbf{d}_{1}^{d-1}$ in a nonzero manner
are $E_1^{*}A_0E_1^{*}$, $E_1^{*}A_{d-1}E_1^{*}, \dots,
E_1^{*}A_{2}E_1^{*}$ in the following manner:
\[E_1^{*}A_0E_1^{*}\mathbf{d}_{1}^{d-1} = \mathbf{d}_{1}^{d-1},\quad
E_1^{*}A_{d-j}E_1^{*}\mathbf{d}_{1}^{d-1} =
2^{j}\mathbf{d}_{1}^{d-1} \mbox{ for } 1 \leq j \leq d-3,\] and
\[E_1^{*}A_{2}E_1^{*}\mathbf{d}_{1}^{d-1}=-
2^{d-2}\mathbf{d}_{1}^{d-1}.\]
\end{enumerate}
It is clear from the above cases that each of $W_{\mathbf{d}_i^l}$
is a $\mathcal{T}$-module. They are irreducible since they are
vector spaces generated by a single vector. This completes the
proof.
\end{proof}
As a consequence, we have the following decomposition of the
standard module.
\begin{theorem} Let $(K_2)^{\wr d}$
be a $d$-class association scheme of order $2^d$. For $l \in \{1,2,
\dots , d-1\}$ define set of vectors $\{\mathbf{d}_{i}^{l}\}_{1 \leq
i \leq 2^{d-l}-1}$ by
\[\mathbf{d}_{i}^{l}=\sum_{k=0}^{2^{l-1}-1}\hat{\mathbf{x}}_{i+j+k}-
\sum_{k=2^{l-1}}^{2^l-1}\hat{\mathbf{x}}_{i+j+k}.\] For each $i$ the
corresponding values of $j$ are successively
\[j=1,\ 1+2^l-1,\ 1+2(2^l-1),\ 1+3(2^l-1),\ \dots,\
1+(2^{d-l}-2)(2^l-1).\] Let $W_{\mathbf{d}_i^l}$ be the linear span
of $\mathbf{d}_{i}^{l}$.

\begin{enumerate}
\item[$(1)$] Let $V$ be the standard module and $\mathcal{P}$ be the primary
module. Then
\[V = \mathcal{P} \oplus \sum W_{\mathbf{v}}\] where $\mathbf{v}$ over all
$\mathbf{d}_{i}^{l}$ defined above.
\item[$(2)$] For $l \in \{1,2,\dots,d-1\}$
\begin{enumerate}
\item[$(a)$] $W_\mathbf{v}$ where $\mathbf{v} \in
    \{\mathbf{d}_{i}^{l}:\ 1 \leq i \leq 2^{d-(l+1)}\}$
    are $\mathcal{T}$-isomorphic.
\item[$(b)$] $W_\mathbf{v}$ where $\mathbf{v} \in
    \{\mathbf{d}_{i}^{l}:\ 2^{d-(l+1)}+1 \leq i \leq
    2^{d-(l+1)}+ 2^{d-(l+2)}\}$ are
    $\mathcal{T}$-isomorphic.
Following above we finally have
\item[$(c)$] $W_{\mathbf{d}_{2^{d-l}-3}^{l}}$ and
    $W_{\mathbf{d}_{2^{d-l}-2}^{l}}$ are
    $\mathcal{T}$-isomorphic.
\item[$(d)$] Rest of $W_\mathbf{v}$ are not
    $\mathcal{T}$-isomorphic.
\end{enumerate}
\end{enumerate}
\end{theorem}
\begin{proof} (1) is straightforward by the construction of
the modules which are mutually orthogonal.\\
(2) The proof is similar when we choose
 $\mathcal{T}$-modules from same groups as described in cases
 (a)-(c). We shall show that $W_{\mathbf{d}_{1}^{l}}$ and
 $W_{\mathbf{d}_{2}^{l}}$ are $\mathcal{T}$-isomorphic.
 Define an isomorphism \[\sigma : W_{\mathbf{d}_{1}^{l}}
 \rightarrow W_{\mathbf{d}_{2}^{l}} \mbox{ by }
\sigma(\mathbf{d}_{1}^{l})=\mathbf{d}_{2}^{l}.\] We need to show
that $(\sigma B - B \sigma)W=0$ for all $B \in \mathcal{T}$. Let us
consider the action of nonzero $E_i^*A_jE_h^*$ on $\sigma$; i.e,
$E_1^{*}A_0E_1^{*}$, $E_1^{*}A_dE_1^{*}$, $E_1^{*}A_{d-1}E_1^{*},
\dots, E_1^{*}A_{d-(l-1)}E_1^{*}$. Now
\[(\sigma E_1^*A_0E_1^*-E_1^*A_0E_1^*\sigma)W_{\mathbf{d}_{1}^{l}}=(\sigma -
E_1^*A_0E_1^* \sigma)W_{\mathbf{d}_{1}^{l}}
%=(\sigma W_{\mathbf{d}_{1}^{l}}) -E_1^*A_0E_1^*(\sigma
%W_{\mathbf{d}_{1}^{l}})
=\mathbf{d}_{2}^{l}-E_1^*A_0E_1^*(W_{\mathbf{d}_{2}^{l}})=0.\]
Similar kind of reasoning shows that for all $E_1^{*}A_0E_1^{*}$,
$E_1^{*}A_dE_1^{*},\ E_1^{*}A_{d-1}E_1^{*},\ \dots,\
E_1^{*}A_{d-(l-1)}E_1^{*}$, $(\sigma E_i^*A_jE_k^*
-E_1^*A_0E_1^*\sigma)W_{\mathbf{d}_{1}^{l}}=0$. Thus,
$W_{\mathbf{d}_{1}^{l}}$ and $W_{\mathbf{d}_{2}^{l}}$ are
$\mathcal{T}$-isomorphic.

Next we shall show that modules selected from different groups are
not $\mathcal{T}$-isomorphic. In particular, let us show that
$W_{\mathbf{d}_{1}^{l}}$ and $W_{\mathbf{d}_{2^{d-(l+1)}+1 }^{l} }$
are not $\mathcal{T}$-isomorphic. Suppose we assume that there
exists an isomorphism
\[\sigma : W_{\mathbf{d}_{1}^{l}} \rightarrow
W_{\mathbf{d}_{2^{d-(l+1)}+1 }^{l} }\mbox{ such that } (\sigma B - B
\sigma)W_{\mathbf{d}_{1}^{l}}=0 \mbox{ for all } B \in
\mathcal{T}.\] Then for $E_1^*A_0E_1^* \in \mathcal{T}$,
\[(E_1^*A_0E_1^*)W_{\mathbf{d}_{1}^{l}}=W_{\mathbf{d}_{1}^{l}}
\mbox{ and } (E_1^*A_0E_1^*)W_{\mathbf{d}_{2^{d-(l+1)}+1 }^{l}
}=0.\] Now,
\[[\sigma(E_1^*A_0E_1^*)-(E_1^*A_0E_1^*)\sigma]W_{\mathbf{d}_{1}^{l}}
=[\sigma-(E_1^*A_0E_1^*)\sigma
]W_{\mathbf{d}_{1}^{l}}=W_{\mathbf{d}_{2^{d-(l+1)}+1 }^{l} } -0 \neq
0\] which is a contradiction to our assumption. Thus,
$W_{\mathbf{d}_{1}^{l}}$ and $W_{\mathbf{d}_{2^{d-(l+1)}+1 }^{l} }$
are not $\mathcal{T}$-isomorphic. The other cases can be proved with
a similar approach.\end{proof}

\begin{theorem} Let $\mathcal{X}=
(K_2)^{\wr d}$ be a $d$-class association scheme of order $2^d$.
\[\mathcal{T} \cong M_{d+1}(\mathbb{C})\oplus M_1(\mathbb{C})^{\oplus
\frac{1}{2} d(d-1)}.\]
\end{theorem}
\begin{proof} It is straightforward. \end{proof}

\section{Concluding Remarks}

There is further work that is needed on the theme related to our
work. Here we state a few problems that are of our interest.
\begin{enumerate}
\item In our attempt to describe the Terwilliger algebra of the
$d$-class association scheme $(K_m)^{\wr d}$ our base was the
Terwilliger algebra described by Tomiyama and Yamazaki \cite{TY94}
for a $2$-class association scheme constructed from a strongly
regular graph. Although our $3$-class association scheme was neither
strongly regular nor a $P$-polynomial scheme, we were able to
describe the Terwilliger algebra concrete manner largely because of
the fact that $d$-class association scheme $(K_m)^{\wr d}$ turned
out to be triply regular, and the structure of the $(d-1)$-class was
beautifully embedded in the $d$-class association scheme. We
demonstrated further how we could extend the same method used to
describe the Terwilliger algebra of the $3$-class association scheme
to the $d$-class association scheme $(K_m)^{\wr d}$ for $d>3$.

The general $d$-class wreath product scheme $K_{n_1}\wr K_{n_2} \wr
\cdots \wr K_{n_d}$ is also triply regular, and its Terwilliger
algebra has the same dimension as in the case of the wreath power
$(K_m)^{\wr d}$. However, it seems to be much more involved to
describe the irreducible $\mathcal{T}$-modules for the general
wreath product scheme. The method used for the case $(K_m)^{\wr d}$
does not work with different $n_i$'s. We do not know how to find
vectors that generate the irreducible $\mathcal{T}$-modules.
However, Paul Terwilliger believes that all non-primary modules
still have dimension $1$. If that is the case, then the structure of
the Terwilliger algebra of $K_{n_1}\wr K_{n_2} \wr \cdots \wr
K_{n_d}$ is also the same as that of $(K_m)^{\wr d}$.

\item In a slightly different direction, it will be interesting to
look at some specific schemes obtained by taking the wreath power of
two association schemes, such as the Hamming $H(2, q)$ instead of
$H(1, q)$. We know that the Terwilliger algebra of a Hamming scheme
can be described as symmetric $d$-tensors on the Terwilliger algebra
of $H(1,q)$ \cite{LMP06}, although in general $H(d, q)$ is not
realized as a product of $H(1, q)$. It would be interesting to see
how the Terwilliger algebra changes when we take the wreath power of
a Hamming scheme $H(d, q)$ for an arbitrary $d>1$.

\item There are also other products besides the wreath
    product. It is also an interesting problem to look at
the direct power of $H(1,q)$. We study the wreath power first
because the direct product of two association schemes has a lot more
classes than the wreath product. Namely, the direct product of a
$d$-class association scheme and a $e$-class association scheme is
of class $de+d+e$ while the wreath product becomes $(d+e)$-class
association scheme. So a study of direct power requires a lot more
work than that of wreath power. However, it may be worthy to look at
it now as we know more about the schemes related to $H(1,q)$.

\item In \cite{Ta97}, the irreducible $\mathcal{T}$-modules
    and Terwilliger algebra has been investigated for the
    Doob schemes. The Doob schemes are the association
    schemes obtained by taking the direct product of copies
    of $H(2,4)$ and copies of schemes coming from the
    Shrikhande graph. In this case, the direct product of
    these schemes preserves many properties of the original
    factor schemes. One important property that is remained
    as the same is $P$-polynomial property. In terms of
    graphs, the Hamming $H(2,4)$ and the Shrikhande graphs
    are the only distance-regular graphs whose direct
    product is also distance-regular. This property no longer
    holds for the direct product of other Hamming Schemes.
    Nevertheless, the
    description of the Terwilliger algebra of Doob
    schemes in terms of those of $H(2,4)$ and Shrikhande
    scheme may shed a light in understanding how the
    Terwilliger algebra of the product behaves when we
    study the Terwilliger algebra of the direct product
    of $H(d,q)$s with various $d$.

\item As explained by Eric Egge \cite{Eg00} and introduction of
\cite{Te92} it is possible to define an ``abstract version" of the
Terwilliger algebra using generators and relations. In all cases the
concrete Terwilliger algebra is a homomorphic image of the abstract
Terwilliger algebra, and in some cases they are isomorphic. In the
case of $(K_m)^{\wr d}$ the entire structure of the Terwilliger
algebra is determined by the intersection numbers and Krein
parameters, so it may be easy to see what is going on. Once all the
vanishing intersection numbers and Krein parameters are worked out,
we can obtain the defining relations for the algebra and we no
longer need to consider the combinatorial structure further.
Terwilliger believes that for the association schemes considered in
the current paper, the abstract Terwilliger algebra and the concrete
Terwilliger algebra are isomorphic. It is remained to study the
Terwilliger algebras (basis, irreducible modules, dimension) from
the generators/relations alone for the association schemes
considered in this paper.

\end{enumerate}

\begin{center}{\textbf{Acknowledgement}}\end{center}
Many thanks are due to Paul Terwilliger who provided many valuable
comments and encouragement during the preparation of this paper.


\begin{thebibliography}{99}

\bibitem{BM95} P. Balmaceda and A. Munemasa, The Terwilliger algebra of group
association scheme, \textit{Kyushu J. Math.}, \textbf{49} (1995),
93--102.

\bibitem{BO94} P. Balmaceda and M. Oura, The Terwilliger
algebra of the group association schemes of $A_5$ and $S_5$, Kyushu
J. Math. \textbf{48} (2) (1994), pp.

\bibitem{BI84} E. Bannai and T. Ito, Algebraic
  combinatorics. I. Association schemes. {\em The Benjamin/Cummings
  Publishing Co., Inc., Menlo Park, CA,} (1984).

\bibitem{Bh08} G. Bhattacharyya, Terwilliger Algebras of
Wreath Products of Association Schemes, Ph.D. Dissertation, Iowa
State University, 2008.

\bibitem{BCN89} A. Brouwer, A.M. Cohen, and A. Neumaier,
\textit{Distancce-Regular Gaphs}, Springer Verlag, Berlin, 1989.

\bibitem{CGS78} P. J. Cameron, J. M. Goethals and J. J.
    Seidel(1978), Strongly regular graphs having strongly regular subconstituents.
     \emph{Journal of Algebra}, (55), 257--280.

\bibitem{Ca99} J. S. Caughman IV, The Terwilliger algebra of
bipartite $P-$ and $Q-$polynomial schemes, \textit{Discrete math.}
\textbf{196} (1999) 65--95.

\bibitem{DC07} B. Curtin and I. Daqqa,
 The subconstituent algebra of a Latin Square.
     \emph{preprint}, 2007

\bibitem{DK94} Yu. A. Drozd and V. V. Kirichenko, \emph{Finite
Dimensional Algebras}, Springer-Verlag, Berlin Heidelberg, 1994.

\bibitem{Eg00} E. Egge, A generalization of the Terwilliger
algebra, \emph{J. Algebra}, {\bf 233} (2000), 213--252.

\bibitem{Go02} J. Go, The Terwilliger algebra of the
hypercube, {\it Europ. J. Combin.} {\bf 23}(4) (2002), 399--429.

\bibitem{Go93} C. D. Godsil, \textit{Algebraic Combinatorics},
Chapman and Hall, New York, 1993.

\bibitem{Hi75} D. G. Higman, Coherent
  configurations. I. Ordinary representation theory.
  {\em Geometriae Dedicata} {\bf 4} (1975), 1--32.


\bibitem{Ja95} F. Jaeger. On spin models, triply regular association schemes, and
duality. \emph{Journal of Algebraic Combinatoris}, (\textbf{2}),
103--144.


\bibitem{LMP06} F. Levstein, C. Maldonado, and D. Penazzi,
The Terwilliger algebra of a Hamming scheme $H(d, q)$. {\it Europ.
J. Combin.} {\bf 27} (2006), 1--10.

\bibitem{Mu93}  A.~Munemasa. An application of Terwilliger algebra.
    \emph{Preprint}, 1993


\bibitem{So02} S.~Y. Song. Fusion Relation in Products of Association Schemes.
     \emph{Graphs and Combinatorics}, (\textbf{18}), 655--665.


\bibitem{Ta97} K. Tanabe.
 The irreducible modules of the Terwilliger algebras of Doob schemes.
     \emph{Journal of Algebraic Combinatorics}, (2), 173--195.


\bibitem{Te92} P. Terwilliger.
    The Subconstituent Algebra of an Association Scheme, (Part I).
     \emph{Journal of Algebraic Combinatorics}, (1), 363--388.

\bibitem{Te93} P. Terwilliger.
    The Subconstituent Algebra of an Association Scheme, (Part II);
    The Subconstituent Algebra of an Association Scheme, (Part III).
     \emph{Journal of Algebraic Combinatorics}, (2), 73--103; 177--210.

\bibitem{Te96} P.~Terwilliger. \textit{Algebraic Combinatorics}.
    Course lecture notes at University of Wisconsin, 1996.

\bibitem{TY94} M. Tomiyama and N. Yamazaki, The subconstituent algebra
of a strongly regular graph, {\it Kyushu J. Math.} {\bf 48} (1994),
323--334.

\end{thebibliography}
\end{document}